\documentclass{amsart}
\usepackage{amscd}

\usepackage{amssymb}
\usepackage{amsmath}
\usepackage{graphicx}
\usepackage{amsfonts}

\newtheorem{theorem}{Theorem}[section]
\newtheorem{proposition}{Proposition}[section]

\newtheorem{remark}{Remark}[section]
\newtheorem{definition}{Definition}[section]
\newtheorem{corollary}{Corollary}[section]
\newtheorem{example}{Example}[section]

\newcommand{\Sing}{\operatorname{Sing}}

\newcommand{\Id}{\operatorname{Id}}

\newcommand{\ox}{\otimes}

\begin{document}
\author{Tornike Kadeishvili}
\thanks{This research described in this publication was made possible in
part by Award No. GM1-2083 of the U.S. Civilian Research and
Development Foundation for the Independent States of the Former
Soviet Union (CRDF) and by Award No. 99-00817 of INTAS}
\address{A. Razmadze Mathematical Institute\\
Georgian Academy of Sciences\\
M. Aleksidze st., 1\\
0193 Tbilisi, Georgia} \email{kade@rmi.acnet.ge}
\author{Samson Saneblidze}
\address{A. Razmadze Mathematical Institute\\
Georgian Academy of Sciences\\
M. Aleksidze st., 1\\
0193 Tbilisi, Georgia} \email{sane@rmi.acnet.ge}
\title{A cubical model for a fibration}
\date{}
\date{ }
\maketitle

\begin{abstract}
In the paper the notion of {\it truncating twisting function} from
a simplicial set to a cubical set and the corresponding notion of
twisted Cartesian product of these sets are introduced. The latter
becomes a cubical set. Using this construction together with the
theory of twisted tensor products for homotopy G-algebras a
strictly associative multiplicative model for a fibration is
obtained.
\end{abstract}

\section{Introduction}

In this paper we construct a {\it  cubical set} which models the
total space of a fibration. The normalized cubical chain complex
of this cubical model {\it coincides} (as a chain complex) with
the twisted tensor product of the singular simplicial complex of
the base and the singular cubical complex of the fiber with
respect to a  certain specific twisting cochain which we call
"truncating". Hence the twisted tensor product  may be endowed
with all structures
 which exist on the chain complex
 of a cubical
set including  the Serre diagonal,  Steenrod chain (co)operations
and other (co)chain operations. In this paper we concentrate only
on the {\it strictly coassociative} Serre diagonal (the cubical
analog of the Alexander-Whitney (AW) diagonal, see \cite{Serre}).
The combinatorial analysis of the Serre diagonal allows us to give
explicit formulas for a strictly associative multiplication on the
twisted tensor product in terms of the $\smile_1$-product and
other related cochain operations measuring the deviation of
$\smile_1$ from being a derivation with respect to the $\smile$
product. Using the standard triangulation of cubes we also obtain
a strictly coassociative diagonal  on Brown's twisted tensor
product of the singular simplicial complex of the base and the
singular simplicial complex of the fiber with respect to some
specific twisting cochain.

For a fibration $F\to E\to Y$, E. Brown \cite{Brown} introduced a {\it %
twisted} differential $d_{\phi}$ on the tensor product
$C^{*}(Y)\otimes C^{* }(F)$ such that the homology of the cochain
complex $(C^{* }(Y)\otimes C^{* }(F),d_{\phi})$ is additively
isomorphic to the cohomology $H^{*}(E)$. There are several papers
(see, for example, L. Lambe and J. Stasheff \cite{Lambe} for
references) where various multiplications are introduced on the
twisted tensor product $C^{*}(Y)\otimes_{\phi} C^{*
}(F)=(C^{*}(Y)\otimes C^{* }(F),d_{\phi})$ to describe
$H^{\ast}(E)$ as an algebra as well. But these multiplications are
either not associative or the differential $d_{\phi}$ is
not a derivation except in special cases, for example, for $Y=S^{n}$ \cite%
{Szczarba}.

The difficulties for introducing  such a multiplication rely on
the following facts. Consider the standard simplicial model of a
fibration: let $X $ be a 1-reduced ($X_{0}=\ X_{1}=pt$) simplicial
set, $G$ a simplicial group, $N$ a simplicial $G$-module,
$t:X_{\ast }\rightarrow G_{\ast -1}$ a twisting function, and
$X\times _{t}N$ the corresponding twisted Cartesian product.
Applying chain functor to $t$ we obtain a twisting cochain
$t_{\ast }=C_{\ast }(t):C_{\ast }(X)\rightarrow C_{\ast -1}(G)$
such that there is a {\it contraction} of $C_{\ast }(X\times
_{t}N)$ to $C_{\ast }(X)\otimes
_{\phi }C_{\ast }(N)$ where $\phi =t_{\ast }$. The simplicial structure of $%
X\times _{t}N$ induces the  AW diagonal on $C_{\ast }(X\times
_{t}N)$. The standard procedure, which uses the basic perturbation
lemma, transports the AW  diagonal to the twisted tensor product
$C_{\ast }(X)\otimes _{\phi }C_{\ast }(N).$ But the resulting
(co)multiplication is (co)associative only up to higher homotopies
\cite{Hueb},\cite{Lambe}.

The situation changes radically if we replace a simplicial group
$G$ by a {\it monoidal cubical set} and suitably modify the notion
of a twisting function. This yields a {\it cubical} model of a
fibration which, as a by-product, induces a strictly associative
multiplication on the above tensor product.

Let us give some more details. Let $X$ be a 1-reduced simplicial
set, $Q$ a {\it monoidal cubical set}, and $L$ a cubical
$Q$-module, i.e., $Q$ and $L$ are cubical sets with given
associative cubical maps $Q\times Q\rightarrow Q$ and $Q\times
L\rightarrow L$. We introduce the notion of {\it truncating
twisting function} $\tau :X_{\ast }\rightarrow Q_{\ast -1}$ from a
simplicial set to a monoidal cubical set (the term {\it
truncating} comes from the universal example $\tau _{U}:\Delta
^{n}\rightarrow I^{n-1}$ of
such functions obtained by the standard truncation procedure, see Section %
\ref{trunc} below). Such a twisting function $\tau $ determines the {\it %
twisted Cartesian product} $X\times _{\tau }L$ as {\it a cubical
set}. We remark that the study of twisting functions from cubical
sets to permutahedral sets and the appropriate twisted Cartesian
product is continued in a forthcoming paper \cite{permu}.

We construct a  functor which assigns to a simplicial set $X$ a
monoidal cubical set ${\bf \Omega} X$ and present a truncating
twisting function $\tau _{U}:X\rightarrow {\bf \Omega} X$ which is
universal in the following sense: Given an arbitrary truncating
function $\tau :X_{\ast }\rightarrow Q_{\ast -1,}$ there is a
monoidal cubical map $f_{\tau }:{\bf \Omega} X\rightarrow Q$ such
that $\tau =f_{\tau }\tau _{U}$. The twisted Cartesian product
${\bf P}X=X\times_{\tau
}{\bf \Omega} X$ is a cubical set that depends functorially on $X$. Note that $%
{\bf \Omega} X$ models the loop space $\Omega |X|$ and ${\bf P}X$
models the path fibration on $|X|$.

The normalized cubical chain functor $C_{\ast }^{\Box }$ applied
to the cubical set ${\bf \Omega} X$ produces $C_{\ast }^{\Box
}({\bf \Omega }X),$ and this chain complex {\it coincides} with
Adams' cobar construction $\Omega C_{\ast }(X)$ (equality (i) of
(\ref{cobarcub})); similarly $C_{\ast }^{\Box }({\bf P}X)$ {\it
coincides} with the acyclic cobar construction $\Omega (C_{\ast
}(X);C_{\ast }(X))$ (equality (ii) of (\ref{cobarcub}));
furthermore $\tau _{\ast }=C_{\ast }(\tau ):C_{\ast
}(X)\rightarrow C_{\ast -1}^{\Box }(Q)$ is a twisting cochain and
$C_{\ast }^{\Box }(X\times _{\tau }L)$ {\it coincides} with the
twisted tensor product $C_{\ast }(X)\otimes _{\tau _{\ast
}}C_{\ast }^{\Box }(L)$ (equality (iii) of (\ref{cobarcub})).

The obtained cubical structures of the cobar construction $\Omega
C_*(X)$ and the twisted tensor product $C_*(X)\otimes_{\tau_*}
C^{\Box}_*(L)$ have the following advantage.

The normalized chain complex of a cubical set  admits the {\it
Serre diagonal} (see \cite{Serre} and below (\ref{SD})), which
turns it into a dg coalgebra.  Since the identification   $C_{\ast
}^{\Box }({\bf\Omega} X)=\Omega C_{\ast }(X)$ the cubical
structure of ${\bf \Omega} X$ determines a  {\it strictly
coassociative comultiplication} on the cobar construction $\Omega
C_{\ast }(X)$. Similarly the cubical structure of $X\times _{\tau
}L$ determines a {\it strictly coassociative comultiplication} on
the twisted tensor product $C_{\ast }^{\Box }(X\times _{\tau
}L)=C_{\ast }(X)\otimes _{\tau _{\ast }}C_{\ast }^{\Box }(L).$
Dually, we immediately obtain the desired strictly associative
multiplication on $C^{\ast }(X)\otimes _{\tau ^{\ast }}C_{\Box
}^{\ast }(L)\subset C_{\Box }^{\ast }(X\times _{\tau }L)$ (here we
have equality when the graded sets have finite type).

Also note that the chain operations dual to Steenrod $\smile_i$
operations are defined for cubical sets in \cite{Kade},
\cite{Kade1} and the equality $C_{\ast }^{\Box }({\bf\Omega}
X)=\Omega C_{\ast }(X)$ allows to define these operations on the
cobar construction $\Omega C_{\ast }(X)$; similarly since
$C_{\ast}^{\Box }(X\times _{\tau }L)=C_{\ast }(X)\otimes _{\tau
_{\ast }}C_{\ast }^{\Box }(L)$ it is possible to introduce
Steenrod operations on multiplicative twisted tensor products.


Next we express the resulting comultiplication on $C_{\ast
}(X)\otimes _{\tau _{\ast }}C_{\ast }^{\Box }(L)$ in terms of
certain chain operations of degree $k$:
\begin{equation*}
E^{k,1}:C_{\ast }(X)\rightarrow C_{\ast }(X)^{\otimes k}\otimes
C_{\ast }(X),\ k\geq 0,
\end{equation*}
which give $C_{\ast }(X)$ a {\it homotopy }G{\it -coalgebra
structure} (dual to a {\it homotopy }G-{\it algebra} in the sense
of Gerstenhaber and Voronov \cite{Voronov}). This structure is a
consequence of the Serre diagonal on $C_{\ast }^{\Box
}({\bf\Omega} X)=\Omega C_{\ast }(X)$: The Serre diagonal of
$C_{\ast }^{\Box}({\bf \Omega} X)$ induces the diagonal $\Omega
C_*(X)\to \Omega C_*(X)\otimes \Omega C_*(X)$ being a
multiplicative map, thus it extends  a certain homomorphism
$C_*(X)\to \Omega C_*(X)\otimes \Omega C_*(X)$, which itself
consists of components $ E^{k,t}:C_{\ast }(X)\rightarrow C_{\ast
}(X)^{\otimes k}\otimes C_{\ast}(X)^{\otimes \ell},\ k,\ell\geq 0,
$ with $E^{k,\ell}=0$ for $\ell\geq 2.$ The operation $E^{1,1}$ is
dual to the Steenrod $\smile _{1}$-cochain operation; thus when
$E^{1,1}=0$ a homotopy G-coalgebra specializes to a cocommutative
dg coalgebra (and dually for homotopy G-algebras). We note that
Baues constructed a homotopy G-coalgebra structure on the
normalized chain complex $C_{\ast }^{N}(X)$ in \cite{Baues1},
\cite{Baues2}.


Towards the end of the paper we develop the theory of
multiplicative twisted tensor products for homotopy G-algebras,
which provides a general algebraic framework for our
multiplicative model of a fibration. First, we review the
theory of multiplicative twisted products due to Prout\`{e} (see \cite%
{Proute}): Suppose $C$ is a dg Hopf algebra, $A$ is a {\it
commutative} dg algebra, $\phi :C\rightarrow A$ is a {\it
coprimitive} twisting cochain (referred to as a {\it
multiplicative} cochain below), and $M$ is simultaneously a dg
algebra and a comodule over $C$ with multiplicative $M\to C\otimes
M$. Then the twisted tensor product $A\otimes _{\phi }M$ is a dga
with respect to the standard multiplication on the tensor product
$A\otimes M$ of dga's. Now replace Prout\`{e}'s commutative $A$ by
a homotopy G-algebra $A$. By definition, there is a strictly
associative multiplication on $BA,$ which can be viewed as a
perturbation of the shuffle product and is compatible with the
coproduct. Thus $BA$ is a dg Hopf algebra. We say that a twisting
cochain $\phi :C\rightarrow A$ is {\it multiplicative} if the
induced map $C\rightarrow BA$
is a dg Hopf algebra map. We introduce a twisted associative multiplication $%
\mu _{\phi }$ on $A\otimes _{\phi }M$ in terms of $\phi $ and the
homotopy G-algebra structure of $A$ by the same formulas as in the
case $A=C^{\ast }(X),$ $C=C_{\Box }^{\ast }(Q)$ and $M=C_{\Box
}^{\ast }(L)$; then $\tau ^{\ast }:C_{\Box }^{\ast }(Q)\rightarrow
C^{\ast }(X)$ provides a basic example of a multiplicative
twisting cochain. Thus, the theory outlined above unifies the
general commutative and homotopy commutative theories; in
particular, this unifies the singular and Sullivan-deRham cochain
complexes of topological spaces.

We remark that the idea of using of {\it cubical} cochains of a
structure group and fiber is found in recent results due to N.
Berikashvili, who constructed a multiplicative model with
associative multiplication when the fiber $F$ is the {\it cubical}
version of an Eilenberg-MacLane space (see \cite{Berika1}) and a
multiplicative model $C^{\ast }(Y)\otimes _{\phi }C_{\Box }^{\ast
}(F),\ \phi :C_{\Box }^{\ast }(G)\rightarrow C^{\ast +1}(Y), $
where $C^{\ast }(Y)$ is the singular {\it simplicial} cochain
complex of the base and $C_{\Box }^{\ast }(G)$ and $C_{\Box
}^{\ast }(F)$ are the singular {\it cubical } cochain complexes of
the structure group and the fiber (see \cite{Berika2}); however,
there is no notion  of underlying truncating twisting functions in
general setting as a map form a simplicial set to a cubical one
leading to the cubical model; also it lacks the analysis of the
Serre cubical diagonal generating the cooperations $E^{k,1},$ and,
consequently,  the general algebraic theory of twisted tensor
products of homotopy commutative dg (co)algebras.

Applying our machinery to a fibration $F\rightarrow E\rightarrow
Y$ on a
1-connected space $Y$ and an associated principal $G$-fibration $%
G\rightarrow P\rightarrow Y$ with action $G\times F\rightarrow F$
we obtain the following cubical model (Theorem \ref{cubmodel}):
Let $X={\operatorname{Sing}}^{1}Y\subset {\operatorname{Sing}}Y$
be the Eilenberg 1-subcomplex generated by the singular simplices
that send the 1-skeleton of the standard $n$-simplex $\Delta ^{n}$
to the base point of $Y.$ Let $Q={\operatorname{Sing}}^{I}G$ and
$M={\operatorname{Sing}}^{I}F$ be the singular cubical sets. Then
Adams' map $\omega _{\ast }:\Omega C_{\ast }(Y)=C_{\ast
}({\bf\Omega} X)\rightarrow C_{\ast }^{\Box }(\Omega Y)$ is
realized by a monoidal cubical map $\omega :{\bf\Omega}
X\rightarrow {\operatorname{Sing}}
^{I}\Omega Y$. Composing $\omega $ with the map of monoidal cubical sets ${%
\operatorname{Sing}}^{I}\Omega Y\rightarrow Q$ induced by the
canonical map $\Omega Y\rightarrow G$ of monoids we immediately
obtain a truncating twisting
function $\tau :X\rightarrow Q$. The resulting twisted Cartesian product ${X}%
\times _{\tau }{M}$ provides the required cubical model of $E;$
and there
exists a cubical weak equivalence ${X}\times _{\tau }{M}\rightarrow {{%
Sing}}^{I}E$. Applying the cochain functor we obtain
Berikashvili's multiplicative twisted tensor product in
\cite{Berika2}.

At the end of the paper we use the theory of multiplicative
twisted tensor products for homotopy G-algebras outlined above to
obtain the multiplicative twisted tensor product $C^{\ast
}(Y)\otimes _{{\phi }}C_{N}^{\ast }(F)$, where $C_{N}^{\ast }$
denotes the normalized singular {\it simplicial}
cochains. The twisting cochain $\phi $ here is the composition ${\phi }%
:C_{N}^{\ast }(G)\overset{\varphi }{\longrightarrow }C_{\Box }^{\ast }(G)%
\overset{\tau ^{\ast }}{\longrightarrow }C^{\ast }(Y)$, where
$\varphi $ is a map of dg Hopf algebras defined by the standard
triangulation of cubes (see the proof of \ref{brownttp}). In other
words, we use a special twisting cochain to introduce an
associative multiplication on Brown's model.

As an example we present fibrations with the base  being a
suspension (in this case the homotopy G-algebra structure consists
just of $E_{1,1}=\smile_1$ and all other operations $E_{k,1}$ are
trivial) and for which the formula for the multiplication in the
twisted tensor product has a very simple form. Moreover in this
case we present small multiplicative model being the twisted
tensor product of cohomologies of base and fiber with the
multiplicative structure purely defined by the $\smile$ and
$\smile_1$ operations.

Finally we mention that the geometric realization $|{\bf\Omega}
\operatorname{Sing}^{1}Y|$ of ${\bf \Omega}
{\operatorname{Sing}}^{1}Y$ is homeomorphic to the cellular model
for a
loop space observed by G. Carlsson and R. J. Milgram \cite{CM}. In \cite%
{Baues1}, \cite{Baues2}, H.-J. Baues defined a geometric
coassociative and homotopy cocommutative diagonal on the cobar
construction $\Omega C_{\ast }^{N}(Y)$ of the normalized chains
$C_{\ast }^{N}(Y)$ by means of a certain cellular model for the
loop space (homotopically equivalent to $|{\bf \Omega }{
\operatorname{Sing}}^{1}Y|$) whose cellular chains coincide with
$\Omega C_{\ast
}^{N}(Y);$ consequently, one obtains a homotopy G-coalgebra structure on $%
C_{\ast }^{N}(Y).$ Another modification of Adams' cobar
construction is considered by Y. Felix, S. Halperin and J.-C.
Thomas \cite{FHT}.

We are indebted to the referee for a number of most helpful
comments and for having suggested many  improvements
 of the exposition.

\section{Notation and preliminaries}

Let $R$ be a commutative ring with unit 1. A {\em differential
graded algebra } (dga) is a graded R-module $C=\{C^{i}\},\
i\in{\mathbb Z},$ with an associative multiplication
$\mu:C^{i}\otimes C^{j}\to C^{i+j}$ and a homomorphism ({\em a
differential} ) $d:C^{i} \to C^{i+1}$ with $d^{2}=0$ and
satisfying the Leibniz rule $d\mu=\mu(d\otimes Id+Id\otimes d)$ .
We assume that a  dga has a unit $\eta :R\to C$ such that
$\mu(\eta\otimes Id)=\mu(Id\otimes \eta)=Id$. A non-negatively
graded
dga $C$ is {\em connected} if $C^{0}=R.$ A connected dga $C$ is {\em %
n-reduced} if $C^{i}=0, 1\leq i\leq n.$ A dga is {\em commutative} if $%
\mu=\mu T,$ where $T(x\otimes y)=(-1)^{|x||y|}(y\otimes x).$ In
general, we use  Koszul's  sign commutation rule: Whenever two
symbols $u$ and $v$ are interchanged, affix the sign
$(-1)^{|u||v|}.$

A {\em differential graded coalgebra } (dgc) is a graded R-module $%
C=\{C_{i}\},\ i\in{\mathbb Z},$ with an coassociative comultiplication $%
\Delta:C\rightarrow C\otimes C$ and a homomorphism ({\em a differential} ) $%
d:C_{i}\rightarrow C_{i-1}$ with $d^{2}=0$ and satisfying $\Delta
d=(d\otimes Id+Id\otimes d)\Delta.$ A dgc $C$ is assumed to have a counit $%
\epsilon :C\rightarrow R,\ \ (\epsilon\otimes
Id)\Delta=(Id\otimes\epsilon)\Delta=Id.$ A non-negatively graded
dgc $C$ is
{\em connected} if $C_{0}=R.$ A connected dgc $C$ is {\em n-reduced} if $%
C_{i}=0,1\leq i\leq n.$ A dgc is {\em cocommutative} if
$\Delta=\Delta T.$

A (connected) {\em differential graded Hopf algebra } (dgha)
$(C,\mu ,\Delta )$ is simultaneously a connected dga $(C,\mu )$
and a connected dgc $ (C,\Delta )$ such that $\Delta :C\rightarrow
C\otimes C$ is an algebra map; note that a graded connected Hopf
algebra has a canonical antipode \cite{Milnor}, so that the
antipode is not an issue.

 A dga $M$ is a (left) {\em comodule}
over a dgha $C$ if $\nu :M\rightarrow C\otimes M$ is a dga map.
Let $\left( M^{\prime },\nu ^{\prime }\right) $ and $\left( M,\nu
\right) $ be comodules over $C^{\prime }$ and $C$, respectively,
and let $\varphi :C^{\prime }\rightarrow C$ be a dgc morphism.
A map $\psi :M^{\prime }\rightarrow M$ is a {\em morphism} of comodules if $%
\nu \psi =(\varphi \otimes \psi )\nu ^{\prime }$.

\subsection{Cobar and Bar constructions}

For an $R$-module $M,$ let $T(M)$ be the tensor algebra of $M$, i.e. $%
T(M)=\oplus _{i=0}^{\infty }M^{\otimes i}$. An element
$a_{1}\otimes ...\otimes a_{n}\in M^{\otimes n}$ is denoted by
$[a_{1},...,a_{n}]$. We denote by $s^{-1}M$ the desuspension of
$M$, i.e. $(s^{-1}M)_{i}=M_{i+1}$.

Let $(C,d_{C},\Delta)$ be a 1-reduced dgc. Denote
$\bar{C}=s^{-1}(C_{>0})$. Let $\Delta=Id\otimes1+1\otimes
Id+{\Delta^{\prime}}$. The (reduced) cobar construction $\Omega C$
on $C$ is the tensor algebra $T(\bar{C})$, with differential
$d=d_{1}+d_{2}$ defined for $\bar{c}\in\bar{C}_{>0}$ by
\begin{equation*}
d_{1}[\bar{c}]=-[\overline{d_{C}(c)}]
\end{equation*}
and
\begin{equation*}
d_{2}[\bar{c}]=\sum(-1)^{|c^{\prime}|}[\bar{c^{\prime}}|\bar{c^{\prime\prime}%
}],\ \ \ \text{for}\ \ \ {\Delta^{\prime}}(c)=\sum
c^{\prime}\otimes c^{\prime\prime},
\end{equation*}
extended as a derivation. The acyclic cobar construction
$\Omega(C;C)$ is the twisted tensor product $C\otimes\Omega C$ in
which the tensor
differential is twisted by the universal twisting cochain $%
C\rightarrow\Omega C$ being an inclusion of degree $-1$ (see
below).

Let $(A,d_{A}, \mu)$ be a 1-reduced dga. The (reduced) bar
construction $BA$ on $A$ is the tensor coalgebra $T(\bar A),\ \bar
A= s^{-1}(A_{>0}),$ with differential $d=d_{1} +d_{2} $ given for
$[\bar a_{1}|\dotsb|\bar a_{n}] \in T^{n}(\bar A)$ by
\begin{equation*}
d_{1}[\bar a_{1}|\dotsb|\bar a_{n}]=-\sum_{i=1}^{n}
(-1)^{\varepsilon_{i}}[\bar
a_{1}|\dotsb|\overline{d_{A}(a_{i})}|\dotsb|\bar a_{n}],
\end{equation*}
and
\begin{equation*}
d_{2} [\bar a_{1}|\dotsb|\bar a_{n}]=- \sum_{i=2}^{n}
(-1)^{\varepsilon_{i}}[\bar
a_{1}|\dotsb|\overline{a_{i-1}a_{i}}|\dotsb|\bar a_{n}],
\end{equation*}
where $\varepsilon_{i}=\sum_{j<i}|\bar{a_{j}}|.$ The acyclic bar
construction $B ( A; A)$ is the twisted tensor product $A\otimes
BA$ in
which the tensor differential is twisted by the universal twisting cochain $%
BA\to A$ being a projection of degree $1.$

\subsection{Twisting cochains}

\label{twist}

Let $(C,d,\Delta :C\rightarrow C\otimes C)$ be a dgc, $(A,d,\mu
:A\otimes A\rightarrow A)$ be a dga, and $(M,d,\nu :M\rightarrow
C\otimes M)$ be\ a dg comodule over $C$. A twisting cochain
\cite{Brown} is a homomorphism $\phi :C\rightarrow A$ of degree 1
satisfying Brown's condition
\begin{equation}
d\phi +\phi d=-\phi \smile \phi ,  \label{brown}
\end{equation}%
where $\phi \smile \phi ^{\prime }=\mu _{A}(\phi \otimes \phi
^{\prime })\Delta _{C}$. There are universal twisting cochains
$C\rightarrow \Omega C$ and $BA\rightarrow A$ being the obvious
inclusion and projection, respectively. Let $T(C,A)$ be the set of
all twisting cochains $\phi
:C\rightarrow A$. Three essential consequences of Brown's condition (\ref%
{brown}) are: \smallskip

\begin{enumerate}
\item[(i)] {\it The multiplicative extension }$f_{\phi }:\Omega
C\rightarrow A${\it \ is a dga map, so there is a bijection}
$T(C,A)\leftrightarrow Hom_{dga}(\Omega C,A)$;\smallskip

\item[(ii)] {\it The comultiplicative extension }$g_{\phi }:C\rightarrow BA$%
{\it \ is a dgc map, so there is a bijection
}$T(C,A)\leftrightarrow Hom_{dgc}(C,BA)$;\smallskip

\item[(iii)] {\it The homomorphism }$d_{\phi }=d\otimes
Id+Id\otimes d+\phi \cap -:A\otimes M\rightarrow A\otimes M${\it ,
where }$\phi \cap (m\otimes
a)=(\mu \otimes Id)(Id\otimes \phi \otimes Id)(Id\otimes \nu )(a\otimes m)$%
{\it , is a differential, i.e.} $d_{\phi }d_{\phi }=0$. \smallskip
\end{enumerate}

\noindent The dg $C$-comodule $(A\otimes M,d_{\phi })$ is called a
twisted tensor product and is denoted by $A\otimes _{\phi }M$. The
twisted tensor product is functorial in the following sense: Let
$\eta :A^{\prime }\rightarrow A$ be a dga morphism, $\varphi
:C^{\prime }\rightarrow C$ be a dgc morphism, $\psi :M^{\prime
}\rightarrow M$ be a morphism of comodules and $\phi ^{\prime
}:C^{\prime }\rightarrow A^{\prime }$ be a twisting cochain such
that $\eta \phi ^{\prime }=\phi \varphi $. Then $\eta \otimes \psi
:A^{\prime }\otimes _{\phi ^{\prime }}M^{\prime }\rightarrow
A\otimes _{\phi }M$ is a chain map.

\subsection{Adams' cobar construction}

\label{Adamscobar}

Let $X$ be a 1-reduced simplicial set, i.e.
$X=\{X_{0}=X_{1}=\{\ast \},X_{2},X_{3},\cdots \}$, and let
$\tilde{C}_{\ast }(X)$ be its chain complex in the ordinary sense.
Define the chain complex ${C}_{\ast }(X)$ as the quotient
\begin{equation*}
{C}_{\ast }(X)=\tilde{C}_{\ast }(X)/\tilde{C}_{>0}(\ast ).
\end{equation*}%
Clearly $C_{\ast }(X)$ is a 1-reduced dgc with respect to the AW
diagonal.

Now let $\operatorname{Sing} Y$ be the singular simplicial set of
a based topological space $Y$ and $X=\operatorname{Sing} ^{1}
Y\subset \operatorname{Sing} Y$ be the (Eilenberg) 1-subcomplex
generated by those singular simplices which send the 1-skeleton of
the standard simplex $\Delta^{n},\ n\geq0,$ to the base point
$y\in Y$. Define the dgc $C_{*}(Y)$ as $C_{*}(X)$. Then Adams'
cobar construction $\Omega C_{*}(Y)$ of a space $Y$ is the cobar
construction of the dgc $C_{*}(Y)$.

\subsection{ Cubical sets}

A cubical set is a graded set $Q=\{Q_{n}\}_{n\geq 0}$ with face operators $%
d_{i}^{\epsilon }:Q_{n}\rightarrow Q_{n-1},\,\epsilon
=0,1,\,i=1,2,...,n,$ and degeneracy operators $\eta
_{i}:Q_{n}\rightarrow Q_{n+1},i=1,2,...,n+1,$ satisfying the
following standard cubical identities \cite{Kan}:


\begin{equation}  \label{ident}
\begin{array}{lll}
d_{j}^{\epsilon }d_{i}^{\epsilon' }=d_{i}^{\epsilon'
}d_{j+1}^{\epsilon },\ \ \ \ \ \ \ \ \, i\leq j
\vspace{1mm}$\newline
$ &  &  \\
d_{i}^{\epsilon }\eta_{j}=\left\{
\begin{array}{ccc}
\eta_{j-1}d_{i}^{\epsilon }\  &  & i<j \\
1 &  & i=j \\
\eta_{j}d_{i-1}^{\epsilon } &  & i>j%
\end{array}
\ \ \right. \vspace{1mm}$\newline
$ &  &  \\
\eta_{i}\eta_{j}=\eta_{j+1}\eta_{i},\ \ \ \ \ \ \ \ \,i\leq j. &  &  \\
&  &
\end{array}%
\end{equation}


\noindent For an example, let $Y$ be a space and let $\operatorname{Sing}^{I}Y=\{%
\operatorname{Sing}_{n}^{I}Y\}_{n\geq 0},$ where
$\operatorname{Sing}_{n}^{I}Y$ is the set of all continuous maps
$I^{n}\rightarrow Y.$ Then $\operatorname{Sing}^{I}Y$ is a cubical
set \cite{Massey}.

Given a cubical set $Q$ and an $R$-module $A$, let $(\bar{C}_{\ast
}^{\Box }(Q;A),d)$ denote its chain complex with coefficients in
$A$. The normalized chain complex $(C_{\ast }^{\Box }(Q;A),d)$ of
$Q$ is defined as the quotient
$C_{\ast }^{\Box }(Q;A)=\bar{C}_{\ast }^{\Box }(Q;A)/D_{\ast }(Q),$ where $%
D_{\ast }(Q)$ is the subcomplex of $(\bar{C}_{\ast }^{\Box
}(Q;A),d)$
generated by the degenerate elements of $Q$. For a space $Y,$ we denote $%
C_{\ast }^{\Box }(\operatorname{Sing}^{I}Y;{\mathbb Z})$ by
$C_{\ast }^{\Box }(Y).$ Both $\bar{C}_{\ast }^{\Box }(Q)$ and
$C_{\ast }^{\Box }(Q)$ are dg coalgebras with respect to the {\it
Serre diagonal} determined by the
Cartesian product  decomposition $I^{n}=I\times \cdots \times I$ of the $n$%
-cube \cite{Serre}: For an element $x\in Q_{n}$ the {\it Serre
diagonal} is given by
\begin{equation}
\label{SD} \Delta (x)=\Sigma (-1)^{\epsilon }d_{j_{1}}^{0}\cdots
d_{j_{p}}^{0}(x)\otimes d_{i_{1}}^{1}\cdots d_{i_{q}}^{1}(x),
\end{equation}
where the summation is over all shuffles
$\{i_{1}<...<i_{q},j_{1}<...<j_{p}\} $ of the set $\{1,...,n\}$
and $(-1)^{\epsilon }$ is the shuffle sign.

Let $Q$ and $Q^{\prime }$ be cubical sets. The {\em (tensor)
product} of $Q$
and $Q^{\prime }$ is defined to be%
\begin{equation*}
Q\times Q^{\prime }=\{(Q\times Q^{\prime
})_{n}=\bigcup_{p+q=n}Q_{p}\times Q_{q}^{\prime }\}/\sim
\end{equation*}%
where $(\eta _{p+1}(a),b)\sim (a,\eta _{1}(b)),\,$ $(a,b)\in
Q_{p}\times Q_{q}^{\prime }.$ This product is endowed with the
obvious face and degeneracy operators \cite{Kan}. Define a {\em
monoidal cubical set} to be a cubical set $Q$ with an associative
cubical multiplication $\mu :Q\times Q\rightarrow Q$ for which a
distinguished element $e\in Q_{0}$ is a unit. (Warning: since the
$Q_{i}$'s are not assumed to be monoids, $Q$ is not a cubical
monoid.) Clearly, the (normalized) chain complex $C_{\ast }^{\Box
}(Q;R)$ on a monoidal cubical set $Q$ and the dual cochain complex
$C_{\Box }^{\ast }(Q;R)$ are dg Hopf algebras. Given a graded
monoidal cubical set $Q,
$ a $Q${\em -module} is a cubical set $L$ together with associative action $%
Q\times L\rightarrow L$ with the unit of $Q$ acting as identity.
In this
case, $C_{\Box }^{\ast }(L;R)$ is a dga comodule over the dg Hopf algebra $%
(C_{\Box }^{\ast }(Q;R),d)$.

\section{The cubical loop and path functors}

\subsection{The cubical loop functor}

In this subsection we construct a functor that assigns to a simplicial set $%
X=\{X_{n},\partial _{i},s_{i}\}$ a cubical monoidal set
${\bf\Omega} X$, which plays the role of the {\it loop space of
}$X$. First we construct a cubical monoid $MX$ {\it without
degeneracies, }then enlarge it to ${\bf \Omega} X$ with degeneracy
operators.

Let $\bar{X}=s^{-1}(X_{>0})$ and define $MX$ to be the free graded
monoid
(without unit) generated by $\bar{X}$. We denote elements of $MX$ by $\bar{x}%
_{1}\dotsm \bar{x}_{k}$ for $x_{j}\in X_{m_{j}+1},\ m_{j}\geq
0,1\leq j\leq k.$ The total degree of an element
$\bar{x}_{1}\dotsm \bar{x}_{k}$ is the
sum $m_{(k)}=m_{1}+\dotsb +m_{k},\ m_{j}=|\bar{x}_{j}|,$ and we write $\bar{x%
}_{1}\dotsm \bar{x}_{k}\in (MX)_{m_{(k)}}$. The product of two elements $%
\bar{x}_{1}\dotsm \bar{x}_{k}$ and $\bar{y}_{1}\dotsm
\bar{y}_{\ell }$ is defined by concatenation $\bar{x}_{1}\dotsm
\bar{x}_{k}\bar{y}_{1}\dotsm \bar{y}_{\ell }$ and is subject only
to the associativity relation; there
are no other relations whatsoever among these expressions. The graded set $%
MX $ canonically admits the structure of a cubical set without
degeneracies in the following fashion: Let
\begin{equation*}
\nu _{i}:X_{n}\rightarrow X_{i}\times X_{n-i},\ \ \nu
_{i}(x)=\partial _{i+1}\dotsm \partial _{n}(x)\times \partial
_{0}\dotsm \partial _{i-1}(x),\ \ 0\leq i\leq n,
\end{equation*}%
denote the components of the AW diagonal. A superscript $n$ on a simplex $%
x^{n}\in X_{n}$ denotes its dimension$.$ Then for an $n$-simplex
$x^{n}\in X_{n},\ n>0,$ let
\begin{equation*}
\nu _{i}(x^{n})=((x^{\prime })^{i},\ (x^{\prime \prime })^{n-i})
\in X_{i}\times X_{n-i}.
\end{equation*}%
First define the face operators
$d_{i}^{0},d_{i}^{1}:(MX)_{n-1}\rightarrow
(MX)_{n-2}$ on a (monoidal) generator $\overline{x^{n}}\in (\bar{X})_{n-1}=%
\overline{X_{n}}$ by
\begin{equation*}
\begin{array}{ll}
d_{i}^{0}(\overline{x^{n}})=\overline{(x^{\prime })^{i}}\cdot \overline{%
(x^{\prime \prime })^{n-i}}, & i=1,...,n-1, \\
d_{i}^{1}(\overline{x^{n}})=\overline{\partial _{i}(x^{n})}, & i=1,...,n-1.%
\end{array}%
\end{equation*}%
Thereafter, for any element (word) $\bar{x}_{1}\dotsm \bar{x}_{k}$
let
\begin{equation*}
\begin{array}{ll}
d_i^0 (\bar x_{1}\dotsm \bar x_{k}) = \bar x_{1}\dotsm \overline {%
(x^{\prime}_q)^{j_q}}\cdot \overline
{(x^{\prime\prime}_q)^{m_q+1-j_q}} \dotsm \bar x_{k} ,$\newline
$\vspace{1mm} &  \\
d^1_i(\bar x_{1}\dotsm \bar x_{k})= \bar x_{1}\dotsm
\overline{\partial _{j_q}(x_q)} \dotsm \bar x_{k}, &
\end{array}%
\end{equation*}%
where $m_{(q-1)}<i\leq m_{(q)},\ j_{q}=i-m_{(q-1)},\ 1\leq q\leq
k,\ 1\leq i\leq n-1.$

It is straightforward to check that the defining identities of a
cubical set hold for $d_{i}^{0},\ d_{i}^{1}$. In particular, the
simplicial relations
between the $\partial _{i}$'s imply the cubical relations between the $%
d_{i}^{1}$'s; the associativity relations between the $\nu _{i}$
's imply the cubical relations between the $d_{i}^{0}$'s, and the
commuting relations between the $\partial _{i}$'s and $\nu _{j}$'s
imply the cubical relations between the $d_{i}^{1}$'s and
$d_{j}^{0}$'s. We now enlarge $MX$ by enlarging its generating set
$\bar{X}$ and introduce the desired degeneracy operators.

For an element $x\in X_{n},$ we consider formal expressions $\eta
_{i_{k}}\cdots \eta _{i_{1}}\eta _{i_{0}}(x)$ with $1\leq
i_{j}\leq n+j-1,\,1\leq j\leq k,\ k\geq 0,$ $\eta _{i_{0}}=Id$. We
call such an expression {\it normal} if $i_{1}\leq \cdots \leq
i_{k}$. Note that any such expression can be reduced to this
normal form by applying the defining
identities for a cubical set with degeneracy operators $\eta _{i}$. Let $%
X^{c}$ be the graded set of formal expressions with normal form
\begin{equation*}
X_{n+k}^{c}=\{\eta _{i_{k}}\cdots \eta _{i_{1}}\eta
_{i_{0}}(x)|\,x\in X_{n}\}_{n\geq 0;k\geq 0},
\end{equation*}%
where
\begin{equation*}
i_{1}\leq \cdots \leq i_{k},\,1\leq i_{j}\leq n+j-1,\,1\leq j\leq
k,\,\eta _{i_{0}}=Id,
\end{equation*}%
and let $\bar{X}^{c}=s^{-1}(X_{>0}^{c}).$ Define ${\bf \Omega}
^{\prime \prime }X$ to be the free graded monoid (without unit)
generated by $\bar{X}^{c}$. It is clear that $X\subset X^{c}$
since $\eta _{i_{0}}(x)=x.$ Thus $MX\subset {\bf \Omega} ^{\prime
\prime }X.$

Let ${\bf \Omega} ^{\prime }X$ be the monoid obtained from ${\bf
\Omega} ^{\prime \prime
}X$ by quotienting with respect to the equivalence relation generated by $%
\overline{\eta _{p+1}(x)}\cdot \bar{y}\sim \bar{x}\cdot
\overline{\eta _{1}(y)}$ for $|x|=p+1,$ $x,y\in X\subset X^{c}.$
We have the inclusion of graded monoids $MX\subset {\bf \Omega}
^{\prime }X$. We claim that ${\bf \Omega} ^{\prime }X$ admits the
structure of a cubical set. Face operators on the subset
$MX\subset {\bf \Omega}^{\prime }X$ were already defined. Now
define a degeneracy operator $\eta _{i}:({\bf \Omega} ^{\prime
}X)_{n-1}\rightarrow ({\bf \Omega}
^{\prime }X)_{n}$ on a (monoidal) generator $\overline{x}\in (\overline{X^{c}%
})_{n-1}$ by
\begin{equation*}
\eta _{i}(\overline{x})=\overline{\eta _{i}(x)},
\end{equation*}%
(assuming $\eta _{i}(x)$ is normalized). For any element
$\bar{x}_{1}\dotsm \bar{x}_{k}$ of ${\bf \Omega} ^{\prime }X$
extend the degeneracy operators by
\begin{equation*}
\begin{array}{ll}
\eta_i(\bar x_{1}\dotsm \bar x_{k})= \bar x_{1}\dotsm \eta_{j_q} (\overline {%
x_q})\dotsm \bar x_{k},$\newline
$\vspace{1mm} &  \\
\eta_n(\bar x_{1}\dotsm \bar x_{k})= \bar x_{1}\dotsm
\bar{x}_{m_{k-1}}\cdot \eta_{m_k+1} (\overline {x_k}), &
\end{array}%
\end{equation*}%
where $m_{(q-1)}<i\leq m_{(q)},\ j_{q}=i-m_{(q-1)},\ 1\leq q\leq
k,$ $1\leq i\leq n-1$. Inductively extend the face operators on
degenerate elements in such a way that the defining identities for
a cubical set are satisfied. Then the cubical set $\{{\bf \Omega}
^{\prime }X,d_{i}^{0},d_{i}^{1},\eta _{i}\}$ depends functorially
on $X$.

Now suppose that $X$ is a based simplicial set with base point
$\ast \in X_{0}$, and denote $e=\overline{s_{0}(\ast )}\in
(\bar{X})_{0}$. Let ${\bf \Omega} X$ be the monoid obtained from
${\bf \Omega} ^{\prime }X$ via
\begin{equation*}
{\bf \Omega }X=\left. {\bf \Omega}^{\prime }X\right/ \sim \text{
,}
\end{equation*}
where $ea\sim ae\sim a,\ $for $a\in {\bf \Omega} ^{\prime }X,$ and $\eta _{n}(\bar{%
x})\sim \overline{s_{n}(x)}$ for $x\in X_{n},$ $n>0.$ Obviously
$({\bf \Omega} X,d_{i}^{0},d_{i}^{1},\eta _{i})$ is a (unital)
monoidal cubical set. Note that although the underlying monoidal
structure of ${\bf \Omega} X$ is not free; all relations involve
degenerate elements.

\begin{remark}
In the definition of the face operators $d_{i}^{0},d_{i}^{1}$ of
${\bf \Omega} X$ for an $n$-simplex of $X_{n}$ the first and last
face operators $\partial
_{0}$ and $\partial _{n}$ of $X$ are not used directly. If, in particular, $%
X $ is a 1-reduced simplicial set (i.e. $X_{0}=X_{1}=\{\ast \}$),
we have
the following identities:%
\begin{equation*}
\begin{array}{ll}
d^0_1(\overline {x^n}) =\overline {(x')^1}\cdot
\overline{(x'')^{n-1}}=e\cdot \overline {(x'')^{n-1}}=
\overline{(x'')^{n-1}}=\overline{\partial_0 (x^n)}, $\newline$\vspace{1mm}\\
d^0_{n-1}(\overline {x^n}) =  \overline {(x')^{n-1}}\cdot
\overline {(x'')^1} =\overline {(x')^{n-1}}\cdot e =
\overline{(x')^{n-1}} = \overline {\partial _n (x^n)}, & x^n\in
X_n.
\end{array}
\end{equation*}%
Thus, all face operators $\partial _{i}$ of $X$ participate in the
definition of ${\bf \Omega} X$ in this case.
\end{remark}

\begin{remark}
The degeneracies of ${\bf \Omega} X$ are formal; we do not use
degeneracies of $X$ except for the last one $s_{n}$. This is
justified by the geometrical fact that in the path fibration, a
degenerate singular n-simplex in the base lifts to a singular
$(n-1)$-cube of the fiber which need not be degenerate (cf. the
proof of Theorem \ref{cobar}).
\end{remark}

It is convenient to verify the cubical relations by the following
combinatorics of the standard cube (compare, \cite{Baues3}).
Motivated by the combinatorial description of the standard $(n+1)$-simplex $%
\Delta ^{n+1},$ we denote the set $\{0,1,...,n+1\}$ by
$[0,1,...,n+1]$ and assign this to the whole $I^{n}.$

\begin{proposition}
\label{cubi} Let
\begin{equation*}
\begin{array}{ll}
d_{i}^{0}\leftrightarrow x_{1},...,x_{i-1},0,x_{i+1},...,x_{n}, &
i=1,...,n
\\
d_{i}^{0}\leftrightarrow x_{1},...,x_{i-1},1,x_{i+1},...,x_{n}, & i=1,...,n%
\end{array}%
\end{equation*}
denote the face operators of the standard cube $I^n$ in Euclidean
coordinates. Then the action of the face operators on
$[0,1,...,n+1]$ by
\begin{equation*}
\begin{array}{lll}
\lbrack 0,1,...,n+1]\overset{d_{i}^{0}}{\longrightarrow} [0,1,...,i][i,...,n+1], & i=1,...,n \\
\lbrack 0,1,...,n+1]\overset{d_{i}^{1}} {\longrightarrow}
 [0,1,...,\hat{i},...,n+1], & i=1,...,n,
\end{array}
\end{equation*}
agrees with the cubical identities.
\begin{proof}
It is straightforward.
\end{proof}
\end{proposition}
 In general,  any q-dimensional
face $a$ of  $I^{n}$ is expressed as
\begin{equation*}
\begin{array}{c}
a=\lbrack 0,i_{1},...,i_{k_{1}}][i_{k_{1}},...,i_{k_{2}}]
[i_{k_{2}},...,i_{k_{3}}]...[i_{k_{p-1}},...,i_{k_{p}},n+1],
\\
0<i_{1}<\ldots <i_{k-{p}}<n+1,\ q=k_{p}-p+1,
\end{array}%
\end{equation*}
in the above  combinatorics; while a cubical degeneracy operator
\begin{equation*}
\eta _{i}\leftrightarrow x_{1},...,x_{i-1},x_{i+1},...,x_{n}
\end{equation*}%
is thought of as adding a formal element $\ast $ to the set
$[0,1,...,n+1]$ at the $(i+1)^{st}$ place:
\begin{equation*}
\eta _{i}[0,1,...,n+1]=[0,1,...,i-1,\ast ,i,...,n+1]
\end{equation*}
with the convention that $[0,1,...,i-1,\ast ][\ast
,i,...,n+1]=[0,1,...,n+1]$ guarantees the equality $d_{i}^{0}\eta
_{i}={Id}=d_{i}^{1}\eta _{i}$.

\subsection{The cubical path functor}

Here we assign to a simplicial set $X$ a cubical set ${\bf P}X$
which plays the role of the {\it path space of} $X$. In some
sense, ${\bf P}X$ will be a {\it twisted Cartesian product} of a
simplicial set $X$ and the monoidal cubical set ${\bf \Omega} X$.

First we define the cubical set ${\bf P}^{\prime }X$ as follows.
Ignoring underlying structure for the moment, consider the
Cartesian product
\begin{equation*}
{X}^{c}{\times }{\bf \Omega} ^{\prime }X=\{({X}^{c}\times {\bf
\Omega} ^{\prime }X)_{n}=\bigcup_{p+q=n}{X}_{p}^{c}\times ({\bf
\Omega} ^{\prime }X)_{q}\}
\end{equation*}%
of the graded sets ${X}^{c}$ and ${\bf \Omega} ^{\prime }X$. Let
\begin{equation*}
{X}^{c}\widetilde{\times }{\bf \Omega} ^{\prime }X={X}^{c}\times
{\bf \Omega} ^{\prime }X/\sim ,
\end{equation*}%
where $(\eta _{p+1}(x),y)\sim (x,\eta _{1}(y)),\ (x,y)\in
{X}_{p}^{c}\times
({\bf \Omega} ^{\prime }X)_{q}$. Introduce operators $d_{i}^{0},d_{i}^{1}$ and $%
\eta _{i}$ on ${X}^{c}\widetilde{\times }{\bf \Omega} ^{\prime }X$
as follows. For
an element $(x,y)\in {\ X}_{p}\times ({\bf \Omega} ^{\prime }X)_{q}\subset {X}%
_{p}^{c}\times ({\bf \Omega} ^{\prime }X)_{q},\ p+q=n,$ let
$$
\begin{array}{ll}
d^0_i  ( x,y) =\begin{cases}
 ((x')^{i-1},\overline {(x'')^{p+1-i}}\cdot y), &
  1\leq i \leq p,  \\

(x,d^0_{i-p}(y)), & p<i\leq n  ,

\end{cases}$\newline$\vspace{0.1in}\\

 d^1_i  ( x,y) =\begin{cases}
(\partial _{i-1}(x),y), &  \ \ \ \ \ \ \ \ \ \ \ \ \ \ 1\leq i\leq p,  \\

(  x, d^1_{i-p}(y)),
 &  \ \ \ \ \ \ \ \ \ \ \  \ \ \ p<i\leq n ,
\end{cases}\\
\end{array}
$$

$$
\begin{array}{lll}
\eta_i(x , y)=(\eta_i(x), y), & 1\leq i \leq
p,$\newline$\vspace{1mm}\\
 \eta_i(x , y)=(x, \eta_{i-p}(y)), & p<i \leq n+1.
\end{array}
 $$
It is easy to check that these face operators satisfy the
canonical cubical identities. The data uniquely extends to the
structure of a cubical set on the whole ${X}^{c}\widetilde{\times
}{\bf \Omega} ^{\prime }X.$ The resulting cubical set is denoted
by ${\bf P}^{\prime }X$; the cubical set ${\bf P}X$ is obtained by
replacing ${\bf \Omega} ^{\prime }X$ by ${\bf \Omega} X$ in the
definition of ${\bf P}^{\prime }X.$ There is the canonical
inclusion of graded sets ${\bf \Omega} X\rightarrow {\bf P}X$
defined by $y\mapsto (\ast ,y),\ \ast \in X_{0},$ and the
canonical projection $\xi :{\bf P}X\rightarrow X$ defined by
$(x,y)\mapsto x.$

 The cubical relations in ${\bf P}^{\prime }X$ can be
verified by means of the following combinatorics of the standard
cube (compare with Proposition \ref{cubi}). The top dimensional
cell of $I^{n+1}$ is identified with the set $0,1,...,n+1]$ while
any proper $q$-face $a$ of $I^{n+1}$ is expressed as
\begin{equation*}
\begin{array}{c}
a=j_{1},...,j_{s_{1}}][j_{s_{1}},...,j_{s_{2}}]
[j_{s_{2}},...,j_{s_{3}}]...[j_{s_{t-1}},...,j_{s_{t}},n+1],
\\
0\leq j_{1}<\ldots <j_{s_{t}}<n+1,\ q=s_{t}-t+1.
\end{array}
\end{equation*}
The dimension of the first block $j_{1},...,j_{s_{1}}]$ is $
\dim([j_{1},...,j_{s_{1}}])+1$.

\begin{proposition}
\label{cubi2}
 Let the face  operators $d^{\epsilon},\,\epsilon=0,1,$
act on   a face  $a$ of $I^{n+1}$ as in Proposition \ref{cubi},
but for its first block as
\begin{equation*}
\begin{array}{ll}
 j_{1},...,j_{s_{1}}]\overset{d^{0}_{i}}{\longrightarrow}
j_1,...,j_{i}][j_{i},...,j_{s_1}], & 1\leq i<s_1,
\vspace{1mm}$\newline$\\
j_{1},...,j_{s_{1}}]\overset{d^{1}_{i}}{\longrightarrow}
j_1,...,\widehat{j_i},...,j_{s_1}],& 1\leq i<s_1.
\end{array}
\end{equation*}
Then the relations among $d^{\epsilon}$'s again  agree with the
cubical identities.
\end{proposition}
\begin{proof}
it is straightforward.
\end{proof}
The canonical cellular map $\psi :I^{n+1}\rightarrow \Delta ^{n+1}$ (\cite%
{Serre}) is combinatorially defined by
\begin{equation*}
j_{1},...,j_{s_{1}}][j_{s_{1}},...,j_{s_{2}}][j_{s_{2}},...,j_{s_{3}}]...[j_{s_{t-1}},...,j_{s_{t}}]\rightarrow
{j_{1},...,j_{s_{1}}}
\end{equation*}%
(see Figure 1). In particular the face $0][0,1,...,n+1]$ of $I^{n+1}$, i.e. $%
d_{1}^{0},$ goes to the minimal vertex (the base point) $0\in
\Delta ^{n+1}.$ \vspace{5mm}

\unitlength=1.00mm \special{em:linewidth 0.4pt}
\linethickness{0.4pt}
\begin{picture}(104.67,86.00)
\put(67.67,85.34){\line(1,0){25.00}}
\put(92.67,85.34){\line(-1,-1){13.67}}
\put(79.00,71.67){\line(-1,0){24.33}}
\put(54.67,71.67){\line(1,1){13.67}}
\put(54.67,71.67){\line(0,-1){22.33}}
\put(54.67,49.34){\line(1,0){24.67}}
\put(79.34,49.34){\line(1,1){13.33}}
\put(92.67,62.67){\line(0,1){22.67}}
\put(79.00,71.67){\line(0,-1){22.33}}
\put(69.00,85.34){\line(0,-1){12.67}}
\put(69.00,69.67){\line(0,-1){7.33}}
\put(54.67,49.34){\line(1,1){7.33}}
\put(63.00,57.34){\line(1,1){6.00}}
\put(69.00,63.00){\line(1,0){9.00}}
\put(80.67,63.00){\line(1,0){12.00}}
\put(15.67,71.34){\line(0,-1){22.00}}
\put(15.67,49.34){\line(1,1){13.67}}
\put(29.34,63.00){\line(0,1){22.33}}
\put(54.67,13.34){\line(3,2){19.00}}
\put(73.67,26.00){\line(3,-2){19.00}}
\put(92.67,13.34){\line(-3,-2){19.00}}
\put(73.67,0.67){\line(-3,2){19.00}}
\put(54.67,13.34){\line(1,0){16.67}}
\put(73.67,25.67){\line(0,-1){25.00}}
\put(75.34,13.34){\line(1,0){16.67}}
\put(15.67,50.00){\circle*{1.33}}
\put(15.67,71.34){\circle*{1.33}}
\put(29.34,62.67){\circle*{1.33}}
\put(54.67,71.67){\circle*{1.33}}
\put(68.67,85.00){\circle*{1.33}}
\put(92.67,85.34){\circle*{1.33}}
\put(79.00,71.67){\circle*{1.49}}
\put(69.00,63.00){\circle*{1.33}}
\put(92.67,63.00){\circle*{1.33}}
\put(79.00,49.34){\circle*{1.33}}
\put(73.67,25.34){\circle*{1.33}}
\put(54.67,13.34){\circle*{1.33}}
\put(92.34,13.34){\circle*{1.33}} \put(73.67,1.34){\circle*{1.33}}
\put(29.34,85.34){\circle*{1.33}}
\put(15.67,71.67){\line(1,1){13.67}}
\put(54.67,49.34){\circle*{1.33}}
\put(33.67,67.67){\vector(1,0){15.67}}
\put(66.67,44.34){\vector(0,-1){14.00}}
\put(58.00,23.67){\vector(-1,1){29.00}}
\put(43.00,43.00){\makebox(0,0)[cc]{$\tau$}}
\put(69.67,37.00){\makebox(0,0)[cc]{$\psi$}}
\put(50.34,13.00){\makebox(0,0)[cc]{$0$}}
\put(75.67,28.00){\makebox(0,0)[cc]{$3$}}
\put(96.00,13.34){\makebox(0,0)[cc]{$2$}}
\put(73.67,-2.66){\makebox(0,0)[cc]{$1$}}
\put(15.67,46.34){\makebox(0,0)[cc]{$0$}}
\put(15.67,75.67){\makebox(0,0)[cc]{$0$}}
\put(32.00,85.67){\makebox(0,0)[cc]{$0$}}
\put(32.00,62.67){\makebox(0,0)[cc]{$0$}}
\put(54.67,46.34){\makebox(0,0)[cc]{$0$}}
\put(79.00,46.34){\makebox(0,0)[cc]{$1$}}
\put(65.34,63.00){\makebox(0,0)[cc]{$0$}}
\put(65.00,85.00){\makebox(0,0)[cc]{$0$}}
\put(96.00,63.00){\makebox(0,0)[cc]{$2$}}
\put(82.34,71.67){\makebox(0,0)[cc]{$1$}}
\put(52.00,71.67){\makebox(0,0)[cc]{$0$}}
\put(104.34,5.34){\makebox(0,0)[cc]{$0123$}}
\put(104.67,69.67){\makebox(0,0)[cc]{$0123]$}}
\put(4.00,65.34){\makebox(0,0)[cc]{$[0123]$}}
\put(95.67,85.33){\makebox(0,0)[cc]{$3$}}
\end{picture}

\vspace{5mm}

\begin{center}
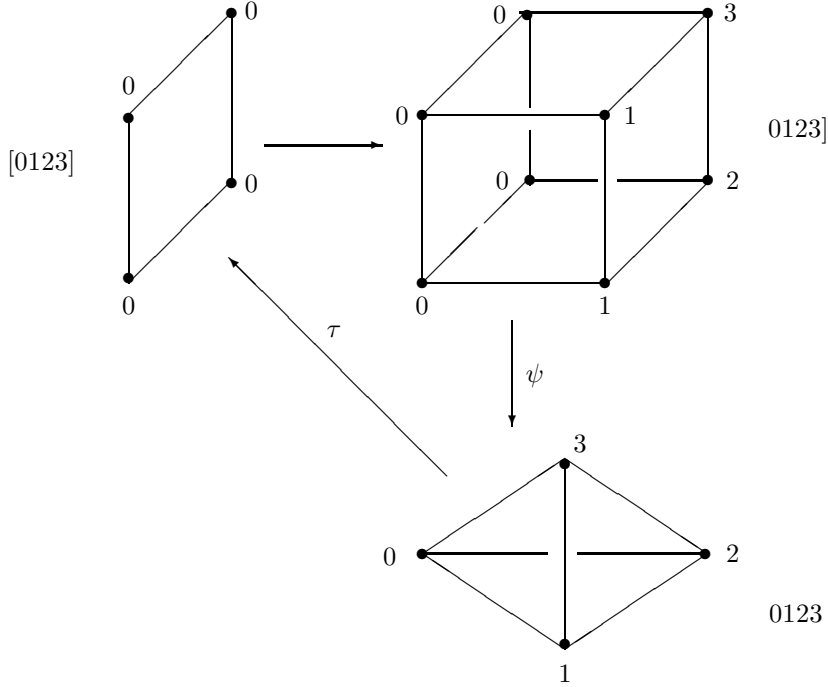
Figure 1: The universal truncating twisting function $\tau.$
\end{center}

\noindent The map $\psi $ can be thought of as a combinatorial
model of the projection ${\bf P}X\overset{\xi }{\longrightarrow
}X$.

\section{Truncating twisting functions and twisted Cartesian products}

\label{trunc}

There is the classical notion of a {\it twisting function} $\tau
:X\rightarrow G$ from a simplicial set to a simplicial group. Such
$\tau $ defines a {\it twisted Cartesian product }for a simplicial
$G$-module $M$ as a simplicial set $X\times _{\tau }M.$ In this
section we introduce the notion of a twisting function between
graded sets in which the domain and the target have face and
degeneracy operators of different types; moreover, the group
structure on each homogeneous component of the target is replaced
by a graded monoidal structure reflecting the standard Cartesian
product of cubes. Namely, we define a {\it truncating twisting
function} $\tau :X\rightarrow Q$ from a simplicial set $X$ to a
monoidal cubical set $Q$. For a cubical $Q$-module with action
$Q\times L\rightarrow L,$ such $\tau $
defines a {\it twisted Cartesian product} $X\times _{\tau }L$ as a {\it %
cubical set}.

These notions are motivated by the cubical set ${\bf P}X,$ which
can be
viewed as a twisted Cartesian product determined by the canonical inclusion $%
\tau :X\rightarrow {\bf \Omega} X,$\ $x\mapsto \bar{x}$ of degree
$-1,$ referred to as the {\em universal truncating twisting
function.}

\begin{definition}
Let $X$ be a 1-reduced simplicial set and $Q$ be a monoidal
cubical set. A sequence of functions $\tau =\{\tau
_{n}:X_{n}\rightarrow Q_{n-1}\}_{n\geq 1} $ of degree $-1$ is
called a truncating twisting function if it satisfies:
\begin{equation*}
\begin{array}{lll}
\tau(x)\ \ \ =e, &  & x\in X_1,$\newline
$\vspace{1mm} \\
d^0_i \tau (x)= \tau \partial_{i+1}\cdots \partial_{n}(x)\cdot
\tau
\partial_{0}\cdots \partial_{i-1}(x), & i=1,...,n-1, & x\in X_n,\, n\geq 1 $%
\newline
$\vspace{1mm} \\
d^1_i\tau (x)=\tau\partial_{i}(x), & i=1,...,n-1, & x\in X_n,\, n\geq 1 $%
\newline
$\vspace{1mm} \\
\eta_{n}\tau(x)=\tau s_n(x), &  & x\in X_{n},\, n\geq 1.%
\end{array}%
\end{equation*}
\end{definition}

\begin{remark}
Note that by definition, a truncating twisting function commutes
only with the last degeneracy operators (compare\thinspace
\cite{Serre}), since this is so for the universal truncating
function.
\end{remark}

The next proposition is an analog of the property (ii) of a
twisting cochain from \ref{twist}.

\begin{proposition}
\label{universal} Let $X$ be a 1-reduced simplicial set and $Q$ be
a monoidal cubical set. A sequence of functions $\tau =\{\tau
_{n}:X_{n}\rightarrow Q_{n-1}\}_{n\geq 1}$ of degree $-1$ is a
truncating twisting function if and only if the monoidal map
$f:{\bf \Omega} X\rightarrow Q$ defined by $f(\bar{x}_{1}\ldots
\bar{x}_{k})=\tau (x_{1})\ldots \tau (x_{k})$ is a map of cubical
sets.
\end{proposition}
\begin{proof}
Since $f$ is completely determined by its restriction to monoidal
generators, use the argument of verification of cubical identities
for a given single generator $\bar{\sigma}$ in ${\Omega} X$ being
equivalent to that of identities of the universal truncating
function $\tau_{U}:\sigma \rightarrow \bar{\sigma}.$
\end{proof}
The following construction is an analog of the property (iii) of a
twisting cochain from \ref{twist}. Given a truncating twisting
function $\tau :X\rightarrow Q$ and a cubical set $L,$ which is a
$Q$-module via $Q\times L\rightarrow L,$ define the corresponding
twisted Cartesian product $X\times _{\tau }L$ by replacing ${\bf
\Omega} X$ with $L$ in the definition of ${\bf P}X$. This gives
the following:

\begin{definition}
\label{tmodule} Let $X$ be a 1-reduced simplicial set, $Q$ be a
monoidal cubical set, and $L$ be a $Q$-module via $Q\times
L\rightarrow L.$ Let $\tau =\{\tau _{n}:X_{n}\rightarrow
Q_{n-1}\}_{n\geq 1}$ be a truncating twisting function. The
twisted Cartesian product $X\times _{\tau }L$ is the graded set
\begin{equation*}
X\times _{\tau }L=X^{c}{\times }L/\sim ,
\end{equation*}%
where $(\eta _{p+1}(x),y)\sim (x,\eta _{1}(y)),\,(x,y)\in
{X}_{p}^{c}\times
L_{q}$, and is endowed with the face $d_{i}^{0},d_{i}^{1}$ and degeneracy $%
\eta _{i}$ operators defined for $(x,y)\in X_{p}\times
L_{q}\subset X_{p}^{c}\times {L}_{q}$ by
\begin{equation*}
\begin{array}{llll}
d^0_i (x,y ) =%
\begin{cases}
(\partial_{1}\cdots \partial_{p}(x), \tau (x)\cdot y) , & i=1, \\
( \partial_{i}\cdots \partial_{p}(x), \tau \partial_{0}\cdots
\partial_{i-2}(x)\cdot y) , & 1<i \leq p, \\
(x,d^0_{i-p}(y)), & p<i\leq n ,%
\end{cases}
$\newline
$\vspace{0.1in} &  &  &  \\
d^1_i ( x,y ) =%
\begin{cases}
(\partial _{i-1}(x),y), & \ \ \ \ \ \ \ \ \ \ \ \ \ \ \ \ \ \ \ \
\ \ \
1\leq i\leq p, \\
( x, d^1_{i-p}(y)), & \ \ \ \ \ \ \ \ \ \ \ \ \ \ \ \ \ \ \ \ \ \
\ p<i\leq
n ,%
\end{cases}
&  &  &  \\
&  &  &
\end{array}%
\end{equation*}
\begin{equation*}
\begin{array}{ll}
\eta_i(x , y)=(\eta_i(x), y), & 1\leq i\leq p, $\newline
$\vspace{1mm} \\
\eta_i(x , y)=(x, \eta_{i-p}(y)), & p<i\leq n+1.%
\end{array}%
\end{equation*}
For any $(x,y)\in X \times _{\tau}{L}$ the operators uniquely
extend to form the cubical set $(X\times _{\tau}L,\
d^0_i,d^1_i,\eta_i) $.
\end{definition}

The geometrical interpretation of $\tau :X\to {\bf \Omega} X$ is
the following: The standard n-simplex (the base) is converted into
the $(n-1)$-cube (the fiber) by the canonical truncation
procedure; this truncation yields the n-cube (the total space) as
well, and the latter is thought of as the "twisted Cartesian
product" of the simplex and the cube (see Fig. 1); so that
projection $\psi$ is a "healing" map. This justifies the name
"truncating twisting function".

\vspace{0.2in}

\begin{example}
Let $M=\{ e_k \}_{k\geq 0}$ be the free graded monoid on a single generator $%
e_1\in M_1$ with trivial cubical set structure and $\tau: X\to M$
the sequence of constant maps $\tau_n:X_n\to M_{n-1},\, n\geq 1.$
Then the
twisted Cartesian product $X\times _{\tau} M$ can be thought of as a {\it %
cubical resolution} of the 1-reduced simplicial set $X.$
\end{example}

The normalized cubical chain functor $C_{\ast }^{\Box }$ applied
to the cubical sets ${\bf \Omega} X,\ {\bf P}X,\ X \times
_{\tau}{L}$ produce dg modules $C_{\ast }^{\Box }({\bf \Omega}X),\
C_{\ast }^{\Box }({\bf P}X),\ C_{\ast }^{\Box }(X \times
_{\tau}{L})$. It is straightforward to check that
\begin{eqnarray}\label{cobarcub}
\begin{array}{lrllll}
 (i)  & C_{\ast }^{\Box }({\bf \Omega} X)&=&\Omega C_{\ast
}(X);\\
\\
 (ii)  & C_{\ast }^{\Box }({\bf P}X)&=&\Omega (C_{\ast}(X);C_{\ast }(X));\\
 \\
 (iii)  & C_{\ast}^{\Box }(X\times _{\tau }L)& = & C_{\ast }(X)\otimes
_{\tau _{\ast }}C_{\ast }^{\Box }(L).
\end{array}
\end{eqnarray}

\section{The cubical model of the path fibration}

Let $Y$ be a topological space. In \cite{Adams}, Adams constructed
a morphism
\begin{equation}
\omega _{\ast }:\Omega C_{\ast }(Y)\rightarrow C_{\ast }^{\Box
}(\Omega Y) \label{adams}
\end{equation}%
of dg algebras that is a weak equivalence for simply connected
$Y$. There are explicit combinatorial interpretations of Adams'
cobar construction, the above map $\omega _{\ast }$, and the
acyclic cobar construction $\Omega (C_{\ast }(Y);C_{\ast }(Y))$ in
terms of cubical sets. Indeed, we have the following theorem
(compare, \cite{Milgram}, \cite{CM}, \cite{Baues1}, \cite
{Baues2}, \cite{FHT}).

\begin{theorem}
\label{cobar} Let $\Omega Y\rightarrow PY\overset{\pi
}{\longrightarrow }Y$ be the Moore path fibration.

(i) There are natural morphisms $\omega ,p,\psi $ such that
\begin{equation}
\begin{CD} \operatorname{Sing} ^I \Omega Y @>>> \operatorname{Sing} ^I P Y
@>\pi _* >> \operatorname{Sing} ^I Y \\ @A \omega AA @A p AA @A {\psi} AA \\
{\bf \Omega} \operatorname{Sing} ^1 Y @>>> {\bf P}
\operatorname{Sing} ^1 Y @>\xi
>> \operatorname{Sing} ^1 Y, \end{CD}  \label{path}
\end{equation}%
$\psi :{\operatorname{Sing}}^{1}Y\rightarrow
{\operatorname{Sing}}^{I}Y$ is a map of graded sets induced by
$\psi :I^{n}\rightarrow \Delta ^{n}$, while $p$ is a morphism of
cubical sets, and $\omega $ a morphism of monoidal cubical sets;
moreover, the cubical maps are homotopy equivalences whenever $Y$
is simply connected.

(ii) The chain complex $C_{\ast }^{\Box }({\bf \Omega}
{\operatorname{Sing}}^{1}Y)$
coincides with the cobar construction $\Omega C_{\ast }(Y)$, see \ref%
{Adamscobar}. Moreover, for a simply connected space, $Y,$ the
Adams weak equivalence (\ref{adams})
$$
\omega _{\ast }:\Omega C_{\ast }(Y)=C_{\ast }^{\Box }(\Omega
{\operatorname{Sing}}^{1}Y) \rightarrow C_{\ast }^{\Box }(\Omega
Y)=C_*(Sing^I\Omega Y)
$$
is induced by the morphism of monoidal cubical sets $\omega $ (and
consequently it preserves all structures which one has  in the
chain complex of a cubical set).

(iii) The chain complex $C_{\ast }^{\Box }({\bf P}
\operatorname{Sing}^{1}Y)
 $  coincides with the acyclic cobar construction\linebreak $\Omega
(C_{\ast }(Y);C_{\ast }(Y)).$
\end{theorem}

\begin{proof}

(i). Morphisms $p$ and $\omega$ are constructed simultaneously by
induction on the dimension of singular  simplices in $\Sing ^1 Y$.
For $i=0,1$ and $(\sigma ,e)\in {\bf P} \Sing ^1 Y$ ,\ $\sigma \in
(\Sing ^1 Y)_i,$ define $p (\sigma ,e)$ as the constant map $I^i
\to P Y$ to the base point $y$,  where $e$ denotes
 the unit of  the monoid ${\bf \Omega}\Sing ^1 Y$ (and of the monoid $\Sing
^I\Omega Y$ as well). Put $\omega (e)=e$. Denote by $P\Sing ^1
Y_{(i,j)}$ the subset in ${\bf P} \Sing ^1Y$ consisting of the
elements $(\sigma, \sigma^{\prime})$ with $|\sigma|\leq i, $  and
$ \sigma^{\prime} \in \Omega \Sing ^1 Y_{(j)}$, a submonoid in
${\bf \Omega}\Sing ^1 Y$ having (monoidal) generators $\bar
\sigma$ with $|\bar \sigma|\leq j$.

Suppose by induction that we have constructed  $p$ and $\omega$ on
${\bf P}  \Sing ^1 Y_{(n-1,n-2)}$ and ${\bf \Omega} \Sing ^1
Y_{(n-2)}$ respectively such that
$$p(\sigma , \sigma^{\prime})= p(\sigma , e)\cdot \omega (\sigma^{\prime})  \ \  \text{and}\ \
\omega (\bar \sigma)= p(d^0_1(\sigma ,e)),
$$
 where the $\cdot$ product  is determined by the action $PY \times \Omega Y\to \Omega Y.$
 Let $\bar I^n \subset I^n$ be the union of the
$(n-1)$-faces $d^{\epsilon}_i(I^n)$  of $I^n$ except the
 $d^0_1(I^n)=(0,x_2,...,x_n)$ and then for a singular simplex
$\sigma: \Delta ^n \to Y$ define the map
$$\bar p: \bar I^n \to   P Y$$
 by $$\bar p|_{d^{\epsilon}_i(I^n)}= p(d^{\epsilon}_{i} (\sigma , e) ),\
\epsilon=0,1,  \ \  \text{and}\ \   i\neq 1\ \  \text{for}\ \
\epsilon=0.   $$ Then the following diagram commutes:
$$
\begin{CD}
\bar I^n @>\bar p>>   P_{\sigma} Y   @> g_{\sigma} >>P   Y \\
@V i VV          @V\pi _{\sigma} VV     @V\pi VV \\
I^n @>\psi >>   \Delta^n    @>\sigma >>   Y.
\end{CD}
$$
Clearly, $i$ is a strong deformation retraction and we define $
p(\sigma, e): I^n \to P_{\sigma }Y  $ as a lift of $ \psi.$ Define
$\omega (\bar \sigma)= p(d^0_1(\sigma ,e)).$ The proof of $p$ and
$\omega$ being homotopy equivalences (after the geometric
realizations) immediately follows, for example, from the
observation that $\xi$ induces a long exact homotopy sequence. The
last statement is a consequence of the following two facts: (1)
$|P\Sing ^1 X|$ is contractible, and (2) the projection $\xi$
induces an isomorphism $\pi_* (|{\bf P} \Sing ^1 Y |, |{\bf
\Omega} \Sing ^1 Y |)\overset{ \xi_*}{\longrightarrow}\pi_*(|\Sing
^1Y|).$

(ii)-(iii). This is straightforward.
\end{proof}

Thus, by passing to chain complexes in diagram (\ref{path}) we
obtain the following comultiplicative model of the path fibration
$\pi$ formed by dgc's.

\begin{corollary}
For the path fibration $\Omega Y\rightarrow PY\overset{\pi
}{\longrightarrow }Y$ there is a comultiplicative model formed by
coassociative dgc's which is natural in $Y$:

\begin{equation}
\begin{CD} C^{\Box}_*( \Omega Y) @>>> C^{\Box}_*( P Y) @>\pi _* >>
C^{\Box}_*(Y) \\ @A \omega_* AA @A p_* AA @A {\psi_*} AA \\ \Omega
C_*( Y) @>>> \Omega(C_*( Y) ;C_*( Y) ) @>\xi_* >> C_*( Y).
\end{CD}
\end{equation}
\end{corollary}

\vspace{0.1in}

\section{Cubical models for fibrations}

Here we prove the main result in this paper. Let $G$ be a
topological group,
$F$ be a $G$-space $G\times F\rightarrow F$, $G\rightarrow P\overset{\pi }{%
\longrightarrow }Y$ be a principal $G$-bundle and $F\rightarrow E\overset{%
\zeta }{\longrightarrow }Y$ be the associated fibration with the
fiber $F$. Let $X={\operatorname{Sing}}^{1}Y$,
$Q={\operatorname{Sing}}^{I}G$ and $L={\operatorname{Sing}}^{I}F.$
The group operation $G\times G\rightarrow G$ induces the structure
of a monoidal cubical set on $Q$ and the action $G\times
F\rightarrow F$ induces a $Q$-module structure $Q\times
L\rightarrow L$ on $L$.

\begin{theorem}
\label{cubmodel} The principal $G$-fibration $G\rightarrow P\overset{\pi }{%
\longrightarrow }Y$ determines a truncating twisting function
$\tau :\operatorname{ Sing}^{1}Y\rightarrow
{\operatorname{Sing}}^{I}G$ such that the twisted Cartesian
product ${\operatorname{Sing}}^{1}Y\times _{\tau
}{\operatorname{Sing}}^{I}F$ models the
total space $E$ of the associated fibration $F\rightarrow E\overset{\zeta }{%
\longrightarrow }Y$, that is there exists a cubical map
\begin{equation*}
{\operatorname{Sing}}^{1}Y\times _{\tau
}{\operatorname{Sing}}^{I}F\rightarrow
{\operatorname{Sing}%
^{I}E}
\end{equation*}%
inducing homology isomorphism.
\end{theorem}

\begin{proof}
Let $\omega:{\bf \Omega} X\to {\Sing }^I\Omega Y$ be the map of
monoidal cubical sets  from Theorem \ref{cobar}. By Proposition
\ref{universal} $\omega$ corresponds to  a truncating twisting
function $\tau':X={\Sing }^1Y
\overset{\tau_{U}}{\longrightarrow}{\bf \Omega} X= {\bf \Omega}
\Sing ^1 Y \overset{\omega} {\longrightarrow}{\Sing }^I\Omega Y.$
 Composing $\tau'$ with the map of monoidal cubical sets ${\Sing
}^I\Omega Y\to {\Sing }^IG=Q $ induced by the canonical map
$\Omega Y\to G$ of monoids we obtain a truncating twisting
function $\tau :X\to Q$. The resulting twisted Cartesian product
${\Sing }^1Y\times_{\tau }{\Sing }^IF$ is a cubical model of $E.$
Indeed, we have the canonical equality
$$X\times_{\tau}L= (X\times_{\tau}Q) \times L/\sim,$$
where $(xg,y)\sim(x,gy).$ Next the argument of the proof of
Theorem \ref{cobar} gives a cubical map $f':X\times_{\tau_U} {\bf
\Omega}X\to \Sing ^I P$ preserving the actions of  ${\bf \Omega}X$
and $Q.$ Hence, this map extends to a cubical map
$f:X\times_{\tau}Q\to \Sing ^I P $ by $f(x,g)=f'(x,e)g.$ The map
$$  f\times Id: (X\times_{\tau}Q) \times  L\to \Sing ^I P\times L\to \Sing ^I (P\times F)$$
induces the map
$$\Sing ^1 Y \times_{\tau} \Sing ^I F\to \Sing ^I E$$
as desired.
\end{proof}


For convenience, assume that $X,Q$ and $L$ are as in the Definition \ref%
{tmodule}. On the chain level a truncating twisting function $\tau
$ induces the twisting cochains $\tau _{\ast }:C_{\ast
}(X)\rightarrow C_{\ast -1}(Q)$ and $\tau ^{\ast }:C^{\ast
}(Q)\rightarrow C^{\ast +1}(X)$ in the standard sense
(\cite{Brown},\cite{Berika3},\cite{Gugenheim}). Recall the
equality of dg modules ((iii) of \ref{cobarcub})
\begin{equation}
C_{\ast }^{\Box }(X\times _{\tau }L)=C_{\ast }(X)\otimes _{\tau
_{\ast }}C_{\ast }^{\Box }(L)  \label{coalgebra}
\end{equation}%
and, consequently, the obvious injection
\begin{equation}
C_{\Box }^{\ast }(X\times _{\tau }L)\supset C^{\ast }(X)\otimes
_{\tau ^{\ast }}C_{\Box }^{\ast }(L)  \label{algebr}
\end{equation}%
of dg modules (which is an equality if the graded sets are of
finite type).

The cubical structure of $X\times _{\tau }L$ induces a dgc structure on $%
C_{\ast }^{\Box }(X\times _{\tau }L)$. Transporting this structure
(the Serre diagonal (\ref{SD})) to the right-hand side of
(\ref{coalgebra}) we obtain a {\it comultiplicative} model
$C_{\ast }(X)\otimes _{\tau _{\ast }}C_{\ast }^{\Box }(L)$ of our
fibration. Dually, $C_{\Box }^{\ast }(X\times _{\tau }L) $ is a
dga, so a dga structure (a multiplication) arises on the
right-hand side of (\ref{algebr}) and we obtain a {\it
multiplicative} model $C^{\ast }(X)\otimes _{\tau ^{\ast }}C_{\Box
}^{\ast }(L)$ of our fibration.

Below we describe these structures (the comultiplication on
$C_{\ast }(X)\otimes _{\tau _{\ast }}C_{\ast }^{\Box }(L)$ and the
multiplication on $ C^{\ast }(X)\otimes _{\tau ^{\ast }}C_{\Box
}^{\ast }(L)$) in terms of certain (co)chain operations that form
a {\it homotopy }G{\it -(co)algebra} structure on the (co)chain
complex of $X$.


\subsection{The canonical homotopy G-algebra structure on $C^*(X)$}

To describe these structures in more detail, we focus on equality
(i) of (\ref{cobarcub})
 \[C_{\ast
}^{\Box }({\bf \Omega} X)=\Omega C_{\ast }(X).
\]
 As before, the cubical structure of ${\bf \Omega} X$ induces a
comultiplication (Serre diagonal) on $C_{\ast }^{\Box }({\bf
\Omega} X)$, thus this structure also appears on the right-hand
side of the above equality, so that the cobar construction $\Omega
C_{\ast }(X)$ becomes a dg Hopf algebra. Such a comultiplication
was defined on the cobar construction $\Omega C_{\ast }^{N}(X)$ of
the normalized complex $C_{\ast }^{N}(X)$ by Baues in
\cite{Baues1}, \cite{Baues2}.

In the combinatorics of Proposition \ref{cubi}, this diagonal is
expressed as
\begin{equation*}
\begin{array}{ll}
\Delta \lbrack 0,1,...,n+1]= & \Sigma (-1)^{\epsilon }\
[0,1,...,j_{1}][j_{1},...,j_{2}][j_{2},...,j_{3}]...[j_{p},...,n+1]\otimes
\\
& [0,j_{1},j_{2},...,j_{p},n+1].
\end{array}%
\end{equation*}%
Note that the summands $[01...n+1]\otimes \lbrack 0,n+1]$ and $
[01][12][23]...[n,n+1]\otimes \lbrack 01...n+1]$ form the
primitive part of the diagonal.

Now regarding the blocks of natural numbers above as faces of the standard $%
(n+1)$-simplex, we obtain Baues' formula for the coproduct $\Delta
:\Omega C_{\ast }(X)\rightarrow \Omega C_{\ast }(X)\otimes \Omega
C_{\ast }(X)$: For
a generator $\sigma \in C_{n+1}(X)\subset \Omega C_{\ast }(X)$ define%
\begin{equation}\label{cohga}
\begin{array}{ll}
\Delta \lbrack \sigma ]= & \Sigma (-1)^{\epsilon }\ [\sigma
(0,1,...,j_{1})|\sigma (j_{1},...,j_{2})|\sigma
(j_{2},...,j_{3})|...|\sigma
(j_{p},...,n+1)]\otimes  \\
& [\sigma (0,j_{1},j_{2},...,j_{p},n+1)],%
\end{array}%
\end{equation}
where $\sigma (i_{1},...,i_{k})$ denotes the suitable face of
$\sigma $. Note that since $X$ is assumed to be 1-reduced, the
image $[\bar{\sigma
(k,k+1)}]$ of each 1-dimensional face $\sigma (k,k+1)$ is the unit in $%
\Omega C_{\ast }(X)$ and hence can be omitted. Note also that the
formula is
highly asymmetric, the left-hand factors of $\Delta \lbrack \sigma ]$ in $%
\Omega C_{\ast }(X)\otimes \Omega C_{\ast }(X)$ have length $\geq
1$ and the right-hand factors have length 1; this is a consequence
of (\ref{SD}) and the structure of $d_{i}^{0},d_{i}^{1}$ from
Proposition \ref{cubi}.

Actually this diagonal consists of {\it components}
\begin{equation*}
E^{k,1}=pr\circ \Delta :C_{\ast }(X)\rightarrow \Omega C_{\ast
}(X)\otimes \Omega C_{\ast }(X)\rightarrow C_{\ast }(X)^{\otimes
k}\otimes C_{\ast }(X),\ k\geq 1,
\end{equation*}%
where $pr$ is the obvious projection. The basic component
$E^{1,1}$ looks like
\begin{equation*}
\begin{array}{lll}
E^{1,1}(\sigma )= & \Sigma _{s,t}(-1)^{\epsilon }\ (\sigma
(0,1)\otimes \sigma (1,2)\otimes ...\otimes \sigma (s-1,s)\otimes
\sigma
(s,s+1,...,t)\otimes \sigma (t,t+1)\otimes  &  \\
& ...\otimes \sigma (n,n+1))\otimes \sigma
(0,1,...,s-1,s,t,t+1,...,n+1)= &
\\
& \Sigma _{s,t}(-1)^{\epsilon }\ \sigma (s,s+1,...,t)\otimes
\sigma (0,1,...,s-1,s,t,t+1,...,n+1), &
\end{array}%
\end{equation*}%
which is a chain operation dual to Steenrod's $\smile
_{1}$-product.

Dualizing the operations $E^{k,1},$ we obtain the sequence of
cochain operations
\begin{equation*}
\{E_{k,1}:C^{\ast }(X)^{\otimes k}\otimes C^{\ast }(X)\rightarrow
C^{\ast }(X)\}_{k\geq 1},
\end{equation*}%
which define a multiplication on the bar construction $BC^{\ast
}(X)\otimes
BC^{\ast}(X)\rightarrow BC^{\ast }(X)$. These cochain operations form a {\it %
homotopy G-algebra structure} on $C^{\ast }(X)$ (see the next
section).

\subsection{The non simply-connected case}

The operations $\left\{ E_{k,1}\right\} $ above are restrictions
of more general cochain operations that arise on $\tilde{C}^{\ast
}(X)$ for a based space $Y$, which is not necessarily 1-connected.
In this case, for $X=\operatorname{ Sing}Y$ we have the operations
\begin{equation*}
\{{E}_{k,1}:\tilde{C}^{\ast }(X)^{\otimes k}\otimes
\tilde{C}^{\ast }(X)\rightarrow \tilde{C}^{\ast }(X)\}_{\ k\geq
0},
\end{equation*}%
given by the following explicit formulas: For $a_{i}\in \tilde{C}%
^{m_{i}}(X),\ m_{i}\geq 2,\ 1\leq i\leq k,$ let
\begin{equation*}
{E}_{k,1}(a_{1},...,a_{k};a_{0})=\sum_{j\geq
k}\tilde{E}_{j,1}(\epsilon ^{1},a_{1},\epsilon ^{1},...,\epsilon
^{1},a_{k},\epsilon ^{1};a_{0}),
\end{equation*}%
where $\epsilon ^{1}\in \tilde{C}^{1}(X)$ is the generator
represented by the constant singular 1-simplex at the base point
$\Delta ^{1}\rightarrow y\in Y$ and the operations
$\tilde{E}_{k,1}$ are defined for $c_{j}\in
\tilde{C}^{m_{j}}(X),\ m_{j}\geq 1,\ 1\leq j\leq k,$ \ $c_{0}\in \tilde{C}%
^{k}(X),$ by
\begin{equation*}
\tilde{E}_{k,1}(c_{1},c_{2},...,c_{k};c_{0})=c\in
\tilde{C}^{n}(X),\ \ \ n=m_{1}+\dotsb +m_{k},
\end{equation*}%
\begin{equation*}
\begin{array}{ll}
c(\sigma)= & (-1)^{\varepsilon}c_1(\partial_{i_1+1}\dotsb
\partial_{n}\sigma) c_2(\partial_0\dotsb \partial_{i_1-1}
\partial_{i_2+1}\dotsb \partial_{n}\sigma)\dotsb $\newline
$\vspace{1mm} \\
& c_k(\partial_{0}\dotsb \partial_{i_{k-1}-1}\sigma)
c_0(\hat\partial_{0}\partial_{1}\hat\partial_{i_1}\dotsb
\hat\partial_{i_{k-1}}\dotsb
\partial_{n-1}\hat\partial_{n}\sigma)$\newline
$\vspace{1mm} \\
& \varepsilon=\sum_{j=1}^{k}(j-1)(m_j-1),%
\end{array}%
\end{equation*}%
where $i_{q}=m_{1}+\dotsb +m_{q},\ 1\leq q\leq k-1,\ \sigma \in
X_{n},$ and where $\tilde{E}_{k,1}(c_{1},c_{2},...,c_{k};c_{0})=0$
otherwise.

\begin{remark}
Though each $\tilde{E}_{k,1},$ and in particular $\tilde{E}_{1,1}$
has only one component, the formula for $k=1$ defines ${E}_{1,1}$
as the Steenrod cochain $\smile_1$-operation without any
restriction on $Y$. This fact evidently indicates a difference
between topological and algebraic interpretation of the operations
$\{E_{k,1}\}_{k\geq 1}$ in terms of 1-reduced algebras (see also
Example \ref{bottsamelson}).
\end{remark}


\subsection{Twisted multiplicative model for a fibration}

Next we further explore the twisted Cartesian product $X\times
_{\tau }L$. To describe the corresponding coproduct and product on
the right-hand sides of (\ref{coalgebra}) and (\ref{algebr})
respectively, it is very convenient to express the Serre diagonal
(\ref{SD}) using the combinatorics of Proposition \ref{cubi2}:
\begin{equation}
\begin{array}{ll}
01...n]\overset{\Delta }{\longrightarrow } & \Sigma (-1)^{\epsilon
}\
0...j_{1}][j_{1}...j_{2}][j_{2}...j_{3}]...[j_{k}...n]\otimes $\newline$\vspace{1mm}\\
& \hat{0},...,\widehat{j_{1}-1},j_{1},\widehat{j_{1}+1},...,\widehat{j_{2}-1}%
,j_{2},...,j_{k},\widehat{j_{k+1}},...,\widehat{n-1},n],%
\end{array}
\label{deltacub}
\end{equation}%
$0\leq j_{1}<\cdots <j_{k}<n,$ where the summands $01...n]\otimes n]$ and $%
0][01][12][23]...[n-1,n]\otimes 01...n]$ form the primitive part
of the diagonal.


Furthermore, the action $Q\times L\rightarrow L$ induces a
comodule structure $\Delta _{L}:C^{\ast }(L)\rightarrow C^{\ast
}(Q)\otimes C^{\ast }(L),$ and it is not hard to see that the
cubical multiplication of (\ref{algebr}) can be expressed by this
comodule structure, the diagonal (\ref{deltacub}), the twisting
cochain $\tau ^{\ast }$, and the operations $\{E_{k,1}\}_{k\geq
1}$ by the following formula: Let $a_{1}\otimes m_{1},\
a_{2}\otimes m_{2}\in C^{\ast }(X)\otimes _{\tau ^{\ast }}C_{\Box
}^{\ast }(L)$ and $\Delta _{L}^{k}:C^{\ast }(L)\rightarrow C^{\ast
}(Q)^{\otimes k}\otimes C^{\ast }(L)$ be the iterated $\Delta _{L}$ with $%
\Delta _{L}^{0}={Id}:C^{\ast }(L)\rightarrow C^{\ast }(L);$ let
$\Delta _{L}^{k}(m_{1})=\sum c^{1}\otimes \ldots \otimes
c^{k}\otimes m_{1}^{k+1}.$ Then
\begin{equation} \label{tformula}
\mu _{\tau ^{\ast }}((a_{1}\otimes m_{1})\otimes (a_{2}\otimes
m_{2}))=\sum_{k\geq 0}(-1)^{|a_{2}||m_{1}^{k+1}|}a_{1}E_{k,1}(\tau
^{\ast }(c^{1}),\dotsc ,\tau ^{\ast }(c^{k});a_{2})\otimes
m_{1}^{k+1}m_{2}.
\end{equation}

\begin{corollary}
\label{twisted} Under the circumstances of Theorem \ref{cubmodel},
the twisted differential $d_{\tau }$ and multiplication $\mu $
turn the tensor product $C^{\ast }(Y)\otimes C_{\Box }^{\ast }(F)$
into a dga $(C^{\ast }(Y)\otimes C_{\Box }^{\ast }(F),d_{\tau
},\mu_{\tau})$ weakly equivalent to the dga $C_{\Box }^{\ast
}(E)$.
\end{corollary}

Such a multiplicative model is constructed in \cite{Berika2}
without explicit formulas for the multiplication.

\begin{corollary}
There exists on the acyclic bar construction $B(C^{\ast
}(Y);C^{\ast}(Y))$ the following strictly associative
multiplication: for $a=a_{0}\otimes \lbrack \bar{a}_{1}|\cdots
|\bar{a}_{n}],\ \ b=b_{0}\otimes \lbrack \bar{b}_{1}|\cdots
|\bar{b}_{m}],\ \ a_{i},b_{j}\in C^{\ast }(Y),\ 0\leq i\leq n,\
0\leq j\leq m,$ let
\begin{equation}
ab=\sum_{k=0}^{n}(-1)^{|b_{0}|(|\bar{a}_{k+1}|+\cdots
+|\bar{a}_{n}|)}a_{0}E_{k,1}(a_{1},\dotsc ,a_{k};b_{0})\otimes
\lbrack \bar{a} _{k+1}|\dotsb |\bar{a}_{n}]\circ \lbrack
\bar{b}_{1}|\dotsb |\bar{b}_{m}]. \label{tbformula}
\end{equation}
\end{corollary}

\begin{proof}
Take $Q=L={\bf \Omega} X$. Then the multiplication
(\ref{tformula}) looks as (\ref{tbformula}).
\end{proof}


\section{Twisted tensor products for homotopy G-algebras}


The notion of homotopy G-(co)algebra naturally generalizes that of
a (co)commutative (co)alge-\linebreak bra. For commutative dga's
there exists the theory of {\it multiplicative} twisted tensor
products. Below we generalize this theory for homotopy G-algebras.
Namely we define a twisted tensor product with both twisted
differential and {\it twisted multiplication } inspired by the
formulas (\ref%
{tformula}) and (\ref{tbformula}) established in the previous
section.


The following definition of homotopy G-algebra (hga) differs from
the definition in \cite{Voronov} only by grading (see also
\cite{GJ}). Let $A$ be a dga and consider the dg module
$({Hom}(BA\otimes BA,A),\nabla )$ with differential $\nabla $. The
$\smile $-product induces a dga structure (the tensor product
$BA\otimes BA$ is a dgc with the standard coalgebra structure).

\begin{definition}
A homotopy G-algebra is a 1-reduced dga $A$ equipped with
multilinear maps
\begin{equation*}
E_{p,q}:A^{\otimes p}\otimes A^{\otimes q}\rightarrow A,\ p,q\geq
0,\ p+q>0,
\end{equation*}%
satisfying the following properties:
\begin{enumerate}
\item[(i)] $E_{p,q}$ is of degree $1-p-q$;\\
\item[(ii)] $E_{p,q}=0$ except $E_{1,0}=id,\ E_{0,1}=id$ and
$E_{k,1},\ k\geq 1$;\\
\item[(iii)] the homomorphism $E:BA\otimes BA\rightarrow A$
defined by
\begin{equation*}
E([\bar{a}_{1}|\dotsb |\bar{a}_{p}]\otimes \lbrack \bar{b}_{1}|\dotsb |\bar{b%
}_{q}])=E_{p,q}(a_{1},...,a_{p};b_{1},...,b_{q})
\end{equation*}
is a twisting cochain in the dga $(Hom(BA\otimes BA,A),\nabla
,\smile )$, i.e. satisfies $\nabla E=E\smile E;$\\
\item[(iv)] the multiplication $\mu _{E}$ is associative, i.e.
$BA$ is a dg Hopf algebra.
\end{enumerate}
\end{definition}

\noindent Condition (iii) implies that the comultiplicative
extension $\mu _{E}:BA\otimes BA\rightarrow BA$ is a chain map;
conditions (iii) and (iv) can be rewritten in terms of the
components $E_{p,q}$ (see \cite{Voronov}). In particular the
operation $E_{1,1}$ satisfies conditions similar to Steenrod's
$\smile _{1}$ product: Condition (iii) gives
\begin{equation}
\begin{array}{c}
dE_{1,1}(a_{1};a_{0})-E_{1,1}(da_{1};a_{0})+(-1)^{|a_{1}|}E_{1,1}(a_{1};da_{0})=
\\
(-1)^{|a_{1}|}a_{1}a_{0}-(-1)^{|a_{1}|(|a_{0}|+1)}a_{0}a_{1},
\end{array}
\label{cup1}
\end{equation}
so it measures the non-commutativity of the product of $A$. Hence,
a homotopy G-algebra with $E_{1,1}=0$ is a commutative dga. We
denote $E_{1,1}(a,b)$ by $a\smile _{1}b$. This notation is also
justified by the other condition that follows from (iii), namely,
\begin{equation}
c\smile _{1}(ab)=(c\smile _{1}a)b+(-1)^{|a|(|c|-1)}a(c\smile
_{1}b). \label{hirsch}
\end{equation}
Thus map $a\smile _{1}-:A\rightarrow A$ is a derivation; when
$A=C^{\ast }(X) $ formula \ref{hirsch} is called the {\it Hirsch
formula}. On the other hand, the map $-\smile _{1}c:A\rightarrow
A$ is a derivation only {\it up to homotopy} with the operation
$E_{2,1}$ serving as a suitable homotopy: This
time condition (iii) gives%
\begin{equation}
\begin{array}{c}
dE_{2,1}(a,b;c)-E_{2,1}(da,b;c)-(-1)^{|a|}E_{2,1}(a,db;c)-(-1)^{|a|+|b|}E_{2,1}(a,b;dc)=
\\
(-1)^{|a|+|b|}(ab)\smile _{1}c-(-1)^{|a|+|b||c|}(a\smile
_{1}c)b-(-1)^{|a|+|b|}a(b\smile _{1}c).%
\end{array}
\label{hirschh}
\end{equation}

The main examples of hga's are: $C^{\ast }(X)$ (see \cite{Baues1}, \cite%
{Baues2},\cite{GJ} and previous section) and the Hochschild
cochain complex of an associative algebra, with the operations
$E_{1,1}$ and $E_{2,1}$ defined by Gerstenhaber in \cite{Gerst}
and the higher operations given in \cite{KadeG},
\cite{GJ},\cite{Voronov}. Another example is the cobar
construction of a dg Hopf algebra \cite{kademeasur}. Note also
that certain algebras (including polynomial algebras) that are
realized as the cohomology of topological spaces also admit a
non-trivial hga structure \cite{saneh}(see also Example
\ref{bottsamelson} below).

The dual notion is that of a homotopy $G$-coalgebra (hgc). For an
hgc $ (C,d,\Delta ,\{E^{p,q}:C\rightarrow C^{\otimes p}\otimes
C^{\otimes q}\})$ the cobar construction $\Omega C$ is a dg Hopf
algebra with a comultiplication induced by $\{E^{p,q}\}$.

\begin{remark}
For a hga $A,$ the operation $E_{2,1}$, besides of (\ref{hirschh}), measures the lack of {\it %
associativity} of $E_{1,1}=\smile _{1}$. In particular, condition
(iv) yields
\begin{equation}
\label{cup1assoc} a\smile _{1}(b\smile _{1}c)-(a\smile
_{1}b)\smile
_{1}c=E_{2,1}(a,b;c)+(-1)^{(|a|+1)(|b|+1)}E_{2,1}(b,a;c),
\end{equation}
which implies that the commutator $[a,b]=a\smile
_{1}b-(-1)^{(|a|+1)(|b|+1)}b\smile _{1}a$ satisfies the Jacobi
identity. In
view of (\ref{cup1}), this commutator induces a Lie bracket of degree -1 on $%
H(A)$. Furthermore, (\ref{hirsch}) and (\ref{hirschh}) imply that $%
[a,-]:H(A)\rightarrow H(A)$ is a derivation, so that $H(A)$ is a
Gerstenhaber algebra \cite{Gerst} (this notion is \underline{not}
a particular case of hga). This structure is generally nontrivial
in the Hochschild cohomology of an associative algebra, but the
existence of a $ \smile _{2}$ product trivializes the induced
Gerstenhaber algebra structure on $H(C^{\ast }(X))=H^{\ast }(X)$.
\end{remark}


\subsection{Multiplicative twisted tensor products}

Let $C$ be a dgc, $A$ a dga and $M$ a dg comodule over $C$.
Brown's twisting
cochain $\phi :C\rightarrow A$ (see \ref{twist}) determines a dga map $%
f_{\phi }:\Omega C\rightarrow A$ (the multiplicative extension of
$\phi $),
a dgc map $g_{\phi }:C\rightarrow BA$ (the comultiplicative extension of $%
\phi $) and the twisted differential $d_{\phi }=d\otimes
Id+Id\otimes d+\phi \cap _{-}:A\otimes M\rightarrow A\otimes M$.
Suppose furthermore, that $C$ is a dg Hopf algebra, $M$ is a dga,
and $M\rightarrow C\otimes M$ is a dga map. In general $d_{\phi }$
is not a derivation with respect to the
multiplication on the tensor product $A\otimes M$. But when $A$ is a {\it %
commutative} dga (in this case $BA$ is a dg Hopf algebra with
respect to the shuffle product $\mu _{sh}$) and $g_{\phi
}:C\rightarrow BA$ is a map of dg Hopf algebras, the twisted
differential $d_{\phi }$ is a derivation with respect to the
standard multiplication of the tensor product $A\otimes C$ and the
twisted tensor product $A\otimes _{\phi }C$ is a dga (see
Prout\`{e} \cite{Proute}).
We shall generalize this phenomenon for a homotopy $G$-algebra
$A,$ in which case $BA$ is again a dg Hopf algebra with respect to
the multiplication $\mu _{E}$.

\begin{definition}
A twisting cochain $\phi :C\rightarrow A$ in $Hom(C,A)$ is
multiplicative if the comultiplicative extension $C\rightarrow BA$
is an algebra map.
\end{definition}

\noindent It is clear that if $\phi :C\rightarrow A$ is a
multiplicative twisting cochain and if $g:B\rightarrow C$ is a map
of dg Hopf algebras then the composition $\phi g:B\rightarrow A$
is again a multiplicative twisting cochain. The canonical
projection $BA\rightarrow A$ provides an example of the universal
multiplicative cochain. For a commutative dga $A$, the
multiplication map $\mu _{E}$ equals $\mu _{sh},$ so Prout\`{e}'s
twisting cochain is multiplicative (see, for example,
\cite{Sane1}). The argument for the proof of formula
(\ref{tformula}) immediately yields the following:

\begin{theorem}\label{ghatwisted}
Let $\phi : C\to A$ be a multiplicative twisting cochain. Then the
tensor product $A\otimes M$ with the twisting differential
$d_{\phi}= d\otimes {Id} +{Id}\otimes d +\phi \cap _{-}$ becomes a
dga $(A\otimes M, d_{\phi},\mu_{\phi })$ with the twisted
multiplication $\mu_{\phi}$ determined by formula
(\ref{tformula}).
\end{theorem}

\begin{remark}
\label{functor} As in \ref{twist}, this construction is functorial
in the following sense: Let $\eta :A^{\prime }\rightarrow A$ be a
{\it strict} morphism of hga's (i.e., $\eta $ is a morphism of
dga's strictly compatible with all $E_{p,q}$'s), $\varphi
:C^{\prime }\rightarrow C$ be a dg Hopf algebra morphism, $\psi
:M^{\prime }\rightarrow M$ be simultaneously a morphism of
comodules and a dga morphism, and $\phi ^{\prime }:C^{\prime
}\rightarrow A^{\prime }$ be a multiplicative twisting cochain such that $%
\eta \phi ^{\prime }=\phi \varphi $. Then
\begin{equation*}
\eta \otimes \psi :(A^{\prime }\otimes M^{\prime },d_{\phi
^{\prime }},\mu _{\phi ^{\prime }})\rightarrow (A\otimes M,d_{\phi
},\mu _{\phi })
\end{equation*}%
is a morphism dga's.
\end{remark}

The above theorem includes the twisted tensor product theory for
commutative algebras (\cite{Proute}).

\begin{corollary}
For a homotopy G-algebra $A$, the acyclic bar-construction
$B(A;A)$, endowed with the twisted multiplication determined by
formula (\ref{tbformula}) acquires a dga structure.
\end{corollary}


\subsection{Brown's model as a dga}

In conclusion, we replace the cubical cochains $C_{\Box }^{\ast }(F)$ and $%
C_{\Box }^{\ast }(G)$ by the normalized {\it simplicial} cochains $%
C_{N}^{\ast }(F)$ and $C_{N}^{\ast }(G)$ in Corollary
\ref{twisted} to introduce an associative multiplication on
Brown's model $C^{\ast
}(Y)\otimes _{\phi }C_{N}^{\ast }(F)$ for a {\it special} twisting cochain $%
\phi $. Specifically, we have:

\begin{corollary}
\label{brownttp} Let $F\rightarrow E\overset{\zeta }{\rightarrow
}Y$ be a fibration as in
Corollary \ref{twisted}. There exists a multiplicative twisting cochain $%
\phi :C_{N}^{\ast }(G)\rightarrow C^{\ast +1}(Y)$ such that the
twisted tensor product $(C^{\ast }(Y)\otimes C_{N}^{\ast
}(F),d_{\phi },\mu _{\phi }) $ with twisted differential $d_{\phi
}$ and twisted multiplication $\mu _{\phi }$ is a dga with
cohomology algebra isomorphic to $H^{\ast }(E)$.
\end{corollary}

\begin{proof} Let us  first mention that
there exists the following standard triangulation of the cub
$I^n$, see for example \cite{FHT}. Each vertex of $I^n$ is a
sequence $(\epsilon _{1},...,\epsilon _{n}),\epsilon _{i}=0,1$.
The set of all $2^{n}$ vertexes is ordered:  $(\epsilon
_{1},...,\epsilon _{n})\leq (\epsilon _{1}^{\prime},...,\epsilon
_{n}^{\prime })$ if $\epsilon _{i}\leq\epsilon _{i}^{\prime} $.
There are $n!$ increasing sequences of maximal length $n+1$ which
start with minimal vertex $(0,...,0)$ and end with maximal
$(1,...,1).$  They form $n!$ $n$-simplices  which triangulate
$I^{n}.$

Let $\varphi :C_{N}^*(G)\rightarrow C_{\Box}^*(G)$ and $\psi
:C_{N}^*(F)\rightarrow C_{\Box}^*(F)$ be  the maps induced by
triangulation of cubes (see, for example, \cite{FHT}), and
$\phi=\tau ^* \varphi:C_N^*(G)\to C_{\Box}^*(G)\to C^*(Y)$. Then
the 4-tuple $\{\eta=Id,\ \varphi,\ \psi,\ \phi \}$ satisfies the
conditions of Remark \ref{functor}, thus
$$
Id\otimes \psi : (C^*(Y)\otimes C_{N}^*(F),d_{\phi},\mu_{\phi})\to
(C^*(Y)\otimes C_{\Box}^*(F),d_{\tau ^* },\mu_{\tau^* })
$$
is a morphism of dga's. A standard spectral sequence argument
shows that this is
 a weak equivalence.
\end{proof}

\subsection{Examples}\label{bottsamelson}
Here we assume that the ground ring $R$ is a field, and all spaces
are path connected. We present examples based on the fact that for
a space being a suspension the corresponding homotopy G-algebra
structure is extremely simple: it consists just of
$E_{1,1}=\smile_1$ and all other operations $E_{k>1,1}$ are
trivial.

  1. The classical Bott-Samelson theorem establishes that
the inclusion $i:X\rightarrow \Omega SX $ induces an algebra
isomorphism
$i_{\ast}:T\tilde{H}_{\ast}(X)\overset{\approx}{\rightarrow}
H_{\ast}(\Omega SX),$ where $SX$ denotes a suspension on a space
$X$. The left hand-side $T\tilde{H}_{\ast}(X)$ is a Hopf algebra
with respect to the comultiplication which extends the one from
$H_{\ast}(X)$ multiplicatively, and the Bott-Samelson map
$i_{\ast}$ is a Hopf algebra isomorphism too. There is the dual
statement for the cohomology as well (cf. Appendix in
\cite{husemoller}).

First we recover the above facts in the following way. Let  $Y$ be
the suspension over a polyhedron $X;$ explicitly, regard $Y$ as
the geometric realization of a quotient simplicial set
$Y=SX/C_{-}X$ where $SX=C_{+}X\cup C_{-}X,$ the union of two cones
over $X$ with the standard simplicial set structure. It is
immediate to check by (\ref{cohga}) that all $E^{k,1}$ for $k\geq
2$ are identically zero, and, moreover, so is the AW diagonal
$\Delta:C_{\ast}(Y)\rightarrow   C_{\ast}(Y) \otimes C_{\ast}(Y)$
in positive degrees as well (cf. \cite{saneh}). Consequently,
since of (\ref{cup1}) and (\ref{cup1assoc})
$E^{1,1}:C_{\ast}(Y)\rightarrow C_{\ast}(Y) \otimes C_{\ast}(Y)$
becomes coassociative chain map of degree 1 and thus it induces a
binary cooperation of degree 1 on the homology denoted by
$Sq^{1,1}:H_{\ast}(Y)\rightarrow H_{\ast}(Y)\otimes H_{\ast}(Y)$.
Notice that both $(C_{\ast}(Y),d,\bar{\Delta}=0, E^{1,1})$ and
$(H_{\ast}(Y),d=0, \bar{\Delta}_{\ast}=0, Sq^{1,1})$ are homotopy
$G$-coalgebras, thus $\Omega C_{\ast}(Y)$ and $\Omega H_{\ast}(Y)$
both are dg Hopf algebras.

The cycle choosing homomorphism $\iota : H_{\ast}(Y)\to
C_{\ast}(Y)$ is a dg coalgebra map in this case. Thus there is a
dg algebra map $\Omega \iota:\Omega H_{\ast}(Y) \to \Omega
C_{\ast}(Y)$ which induces the Bott-Samelson isomorphism of graded
algebras
\begin{equation}
\label{Bott} T \tilde{H}_{\ast}(X)=\Omega H_{\ast}(Y)= H(\Omega
H_{\ast}(Y))
\overset{(\Omega\iota)_{\ast}}{\longrightarrow}H_{\ast}(\Omega
C_{\ast}(Y))=H_{\ast}(\Omega Y).
\end{equation}
To show that (\ref{Bott}) is a Hopf algebra isomorphism, let first
consider  the  diagram
\[
\begin{CD}
C_{\ast}(Y)@>{E^{1,1}}>> C_{\ast}(Y) \otimes  C_{\ast}(Y) \\
@AA{\iota}A  @AA{\iota}\ox{\iota}A    \\
\tilde{H}_{\ast}(Y) @>{Sq^{1,1}}>> \tilde{H}_{\ast}(Y) \otimes
\tilde{H}_{\ast}(Y)\\
@A{\approx}A s A   @A{\approx}A s\otimes s A\\
\tilde{H}_{\ast-1}(X) @>\Delta_{\ast}>> \tilde{H}_{\ast-1}(X)
\otimes \tilde{H}_{\ast-1}(X),
\end{CD}
\]
where $s$ is  the suspension isomorphism;  the upper square is
commutative up to a chain homotopy, while the bottom square is
strict commutative. This implies that $\Omega \iota$ is also a
coalgebra map up to a chain homotopy, consequently (\ref{Bott}) is
a  coalgebra map too.

2. Let
 $\Omega Y\rightarrow PY\overset{\pi}{\rightarrow} Y$
be the Moore path fibration with the base $Y$ which is the
suspension over a polyhedron $X.$ Let $f:Y\rightarrow Z$ be a map,
$ \Omega Y\times \Omega Z\rightarrow \Omega Z$ be the induced
action via the composition
\[\Omega Y\times \Omega Z \overset{\Omega f\times \Id}{\longrightarrow}
\Omega Z\times \Omega Z\rightarrow \Omega Z,\] and  $\Omega
Z\rightarrow E_{f}\overset{\xi}{\rightarrow} Y$ be the associated
fibration; for simplicity assume that
  $Z$ is
the suspension and  simply connected $CW$-complex of finite type,
as well. We present two multiplicative models for the fibration
$\xi$ using the cubical model $Y\times_{\tau} {\bf\Omega} Z$ with
the universal truncating twisting function $\tau=\tau_{U}:
Y\rightarrow  {\bf \Omega} Y$.

Notice that the twisted  differential of the cochain complex
$(C^{\ast}(Y\times_{\tau} {\bf\Omega} Z),d)=(C^{\ast}(Y)\otimes
C^{\ast}({\bf \Omega} Z),d_{\tau^{\#}})=(C^{\ast}(Y)\otimes B
C^{\ast}( Z),d_{\tau^{\#}})$ with universal
$\tau^{\#}:BC^{\ast}(Y)\rightarrow C^{\ast}(Y)$ becomes  the form
$$
\begin{array}{c}
d_{\tau^{\#}}(a\otimes [\bar{m}^1|...|\bar{m}^n]) = da\otimes
[\bar{m}^1|...|\bar{m}^n] +\sum_{k=1}^{n} a\otimes
[\bar{m}^1|...|d\bar{m}^k|...|\bar{m}^n]+
 \\
a\cdot m_1\otimes [\bar{m}^2|...|\bar{m}^n].
\end{array}
$$
Since the simplified  structure of the homotopy G-algebra
$(C^{\ast}(Y),d,\mu=0, E_{1,1})$  formula (\ref{tformula}) becomes
the following form:
\begin{equation}
\label{multc} \mu_{\tau^{\#}}((a_1\otimes m_1)(a_2\otimes
m_2))=a_1a_2\otimes m_1m_2+a_1E_{1,1}(f^{\#}(m^1_1),a_2)\otimes
[\bar{m}^2_1|...|\bar{m}^n_1]\cdot m_2,
\end{equation}
where $f^{\#}:C^{\ast}( Z)\rightarrow C^{\ast}( Y) ,$\,
$a_1,a_2\in C^{\ast}(Y), \,
m_1=[\bar{m}^1_1|...|\bar{m}^n_1],m_2\in B C^{\ast}( Z),\, n\geq
0.$ Note that since the product on ${C}^{>0}(Y)$ is zero, the
twisted part of $\mu_{\tau^{\#}}$ (the second summand) may be
non-zero only for $a_1\in C^{0}(Y).$

So that we  get that $H(C^{\ast}(Y)\otimes B C^{\ast}(
Z),d_{\tau^{\#}},\mu_{\tau^{\#}})$ and $H^*(E_f)$ are isomorphic
as algebras.

On the other hand, let us consider the following multiplicative
twisted tensor product $(H^{\ast}(Y)\otimes H^{\ast}({\bf \Omega}
Z),d_{\tau^{\ast}})=(H^{\ast}(Y)\otimes B H^{\ast}(
Z),d_{\tau^{\ast}})$ with universal
$\tau^{\ast}:BH^{\ast}(Y)\rightarrow H^{\ast}(Y)$. The
differential here is of the form:
$$
d_{\tau^{*}}(a\otimes [\bar{m}^1|...|\bar{m}^n]) = a\cdot
m_1\otimes [\bar{m}^2|...|\bar{m}^n].
$$
Again since the simplified structure of the homotopy G-algebra
$(H^{\ast}(Y),d=0, \mu^{\ast}=0, Sq_{1,1})$ the formula
(\ref{tformula}) becomes the following form:
\begin{equation}
\label{multh} \mu_{\tau^{\ast}}((a_1\otimes m_1)(a_2\otimes
m_2))=a_1a_2\otimes m_1m_2+a_1Sq_{1,1}(f^{\ast}(m^1_1),a_2)\otimes
[\bar{m}^2_1|...|\bar{m}^n_1]\cdot m_2,
\end{equation}
where $ f^{\ast}:H^{\ast}(Z)\rightarrow H^{\ast}(Y),\, a_1,a_2\in
H^{\ast}(Y),\, $ $m_1=[\bar{m}^1_1|...|\bar{m}^n_1],m_2\in B
H^{\ast}( Z),\, n\geq 0.$ Note that since the product on
${H}^{>0}(Y)$ is zero,  the twisted part of $\mu_{\tau^{\ast}}$
(the second summand) may be non-zero only for $a_1\in {H}^{0}(Y).$
Also we remark that for  an element $a\in H^{\ast}(Y),$ one gets
$Sq_{1,1}(a,a)=Sq_1(a),$ the Steenrod square.

We claim that $(H^{\ast}(Y)\otimes B H^{\ast}(Z),d_{\tau^{\ast}})$
is a "small" multiplicative model of the fibration $\xi,$ i.e.,
$H(H^{\ast}(Y)\otimes B H^{\ast}(Z),d_{\tau^{\ast}})$ and
$H^*(E_f)$ are isomorphic as algebras. Indeed, it is
straightforward to calculate  (or using the standard spectral
sequence argument) that additively
$$
\begin{array}{c}
H(C^{\ast}(Y)\otimes B C^{\ast}( Z),d_{\tau^{\#}})\approx
H(H^{\ast}(Y)\otimes B H^{\ast}( Z),d_{\tau^{\ast}})\approx \\
H^{0}(Y)\otimes T_f(H^{\ast}(Z))\bigoplus
H^{\ast}(Y)/Imf^{\ast}\otimes BH^{\ast}(Z) ,
\end{array}
$$
where $T_f(H^{\ast}(Z))=s^{-1}(Ker f^{\ast})+ s^{-1}(Ker
f^{\ast})\otimes s^{-1}H^{\ast}(Z)+\cdots+ s^{-1}(Ker
f^{\ast})\otimes
 (s^{-1}H^{\ast}(Z))^{\otimes n}
+\cdots ,\, n\geq 1.$ Since the explicit formulas (\ref{multc})
and (\ref{multh}) it is easy to calculate that the twisted parts
of $\mu_{\tau^{\#}}$ and $\mu_{\tau^{\ast}}$ annihilate in
homology, thus they induce the same multiplication on
$H^{\ast}(E_f).$ As a byproduct we obtain that the multiplicative
structure of the total space $E_f$ does not depend on a map $f$ in
a sense that if $f^{\ast}=g^{\ast}$ then
$H^{\ast}(E_f)=H^{\ast}(E_g)$ as algebras.  Note also that this
multiplicative structure is purely defined by the $\smile$ and
$\smile_1$ operations.


\vspace{0.2in}


\begin{thebibliography}{99}
\bibitem{Adams} J. F. Adams, On the cobar construction, Proc. Nat. Acad.
Sci. (USA), 42 (1956), 409-412.


\bibitem{Baues1} H.-J. Baues, Geometry of loop spaces and the cobar
construction, Memoires of the AMS, 25 (1980), 1-170.

\bibitem{Baues2} -----------, The double bar and cobar construction,
Compositio Math., 43 (1981), 331-341.

\bibitem{Baues3} -----------, The Cobar Construction as a Hopf Algebra,
Invent. Math., 132 (1998) 467-489 .

\bibitem{Berika1} N. Berikashvili, An algebraic model of fibration with the
fiber $K(\pi ,n)$-space, Georgian Math. J., 3, 1 (1996), 27-48.

\bibitem{Berika2} N. Berikashvili and D. Makalatia, The multiplicative
version of twisted tensor product theorem, Bull. Georg. Acad.
Sci., 154 (1996), 327-329.

\bibitem{Berika3} N. Berikashvili, On the differentials of spectral
sequences (Russian), Proc. Tbilisi Mat. Inst., 51 (1976), 1-105.

\bibitem{Brown} E. Brown, Twisted tensor products, Ann. of Math., 69 (1959),
223-246.

\bibitem{CM} G. Carlsson and R. J. Milgram, Stable homotopy and iterated
loop spaces, Handbook of Algebraic Topology (Edited by I. M.
James), North-Holland (1995), 505-583.

\bibitem{FHT} Y. Felix, S. Halperin and J.-C. Thomas, Adams' cobar
equivalence, Trans. Amer. Math. Soc., 329 (1992), 531-549.

\bibitem{Gerst} M. Gerstenhaber, The cohomology structure of an associative
ring, Ann. of Math., 78 (1963), 267-288.

\bibitem{Voronov} M. Gerstenhaber and A. A. Voronov, Higher Operations on
the Hochschild Complex, Functional Analysis and its Applications,
29 (1995), 1-5.

\bibitem{GJ} E. Getzler and J. D. Jones, Operads, homotopy algebra, and
iterated integrals for double loop spaces, preprint, 1995.

\bibitem{Gugenheim} V.K.A.M. Gugenheim, On the chain complex of a fibration,
Ill. J. Math., 16 (1972), 398-414.

\bibitem{Hueb} J. Huebschmann, The homotopy type $F\Psi^q,$ the complex and
sympletic cases, Cont. Math., 55 (1986), 487-518.

\bibitem{husemoller} D. Husemoller, Fibre bundles, third ed.,
Graduate Texts in Mathematics, no. 20, Springer-Verlag, New-York,
1994.


\bibitem{Kan} D. M. Kan, Abstract homotopy I, Proc. Nat. Acad. Sci. U.S.A.,
41 (1955), 1092-1096.

\bibitem{Kade} T. Kadeishvili, DG Hopf Algebras with Steenrods i-th
coproducts, Bull. Georgian Acad. Sci., 158 (1998), 203-206.

\bibitem{Kade1} -------------, Cochain operations defining Steenrod $%
\smile_i $-products in the bar construction, Georgian Math.
Journal, v. 10 (2003), 115-125.

\bibitem{KadeG}

------------, The $A(\infty )$-algebra structure and cohomology of
Hochschild and Harrison (Russian), Proc. Tbilisi Math. Inst., 91
(1988), 19-27.

\bibitem{kademeasur} ------------, Measuring the noncommutativty of DG
algebras, preprint.

\bibitem{permu} T. Kadeishvili and S. Saneblidze, The twisted Cartesian
model for the double path space fibration, preprint, AT/0210224.


\bibitem{Lambe} L. Lambe and J. Stasheff, Applications of perturbation
theory to iterated fibrations, Manusc. Math., 58 (1987), 363-376.

\bibitem{Massey} W. S. Massey, A basic course in Algebraic topology,
Springer-Verlag, Ney-York, 1991.

\bibitem{Milgram} R. J. Milgram, Iterated loop spaces, Ann. of Math., 84
(1966), 386-403.

\bibitem {Milnor} J. Milnor and J. C. Moore, On the Structure of Hopf Algebras,
   Ann. of Math. 81 (1965), 211-264.


\bibitem{Proute} A. Proute, $A_{\infty}$-structures, Modele minimal de
Bauess-Lemaire des fibrations, preprint.

\bibitem{Sane1} S. Saneblidze, Perturbation and obstruction theories in the
fibre spaces, Proc. A. Razmadze Math. Inst., 111 (1994), 1-106.

\bibitem{saneh} ------------, The Hochschild complex of a space is the
complex of the Hochschild set, preprint.

\bibitem{Serre} J.-P. Serre, Homologie singuliere des espaces fibres,
applications, Ann. Math., 54 (1951), 429-505.

\bibitem{Szczarba} R. H. Szczarba, The homology of twisted cartesian
products, Trans. AMS, 100 (1961), 197-216.
\end{thebibliography}
\end{document}